\documentclass[12, final]{article}
\usepackage{graphicx}
\usepackage[thmmarks]{ntheorem}
\usepackage{amsmath, amssymb}
\usepackage{booktabs}
\newcommand{\bs}[1]{\boldsymbol{#1}}

\usepackage{natbib}

\usepackage{xcolor}
\usepackage{lscape}

\newtheorem{thm}{Theorem}
\newtheorem{lem}{Lemma}

\newtheorem{assn}{Assumption}

\theoremstyle{nonumberplain} \theorembodyfont{\textrm{\text}}
\theoremsymbol{\ensuremath{_\square}}

\begin{document}

\thispagestyle{empty}

\title{\Large \bf A $U$-classifier for high-dimensional data under non-normality}

\author{{\sc M. Rauf Ahmad}$^{1}$ and {\sc Tatjana Pavlenko}$^{2}$\\[0.1cm]
$^1$ {\small \it Department of Statistics, Uppsala University, Uppsala, Sweden}\\
$^2$ {\small \it Department of Mathematics, KTH Royal Institute of Technology,}\\{\small \it Stockholm, Sweden}}

\date{}
\maketitle

\begin{abstract}
A classifier for two or more samples is proposed when the data are high-dimensional and the underlying distributions may be non-normal. The classifier is constructed as a linear combination of two easily computable and interpretable components, the $U$-component and the $P$-component. The $U$-component is a linear combination of $U$-statistics which are averages of bilinear forms of pairwise distinct vectors from two independent samples. The $P$-component is the discriminant score and is a function of the projection of the $U$-component on the observation to be classified. Combined, the two components constitute an inherently bias-adjusted classifier valid for high-dimensional data. The simplicity of the classifier helps conveniently study its properties, including its asymptotic normal limit, and extend it to multi-sample case. The classifier is linear but its linearity does not rest on the assumption of homoscedasticity. Probabilities of misclassification and asymptotic properties of their empirical versions are discussed in detail. Simulation results are used to show the accuracy of the proposed classifier for sample sizes as small as 5 or 7 and any large dimensions. Applications on real data sets are also demonstrated.
\end{abstract}

\vspace{0.5cm}


%%%%%%%%%%%%%%%%%%%%%%%%%%%%%%%%%%%%%%%%%%%%%%%%%%%%%%%%%%%%%%%%%%%%%%%%%%%%%%%%%%%%%%%%%%%%%%%
\section{Introduction}\label{sec:Intro}
%%%%%%%%%%%%%%%%%%%%%%%%%%%%%%%%%%%%%%%%%%%%%%%%%%%%%%%%%%%%%%%%%%%%%%%%%%%%%%%%%%%%%%%%%%%%%%%

A linear classifier for two or more populations is presented when the data are high-dimensional and possibly non-normal. Let
\[
{\bf x}_{ik} = (x_{i1k}, \ldots, x_{ikp})', ~~k = 1, \ldots, n_i,
\]
be $n_i$ independent and identically distributed random vectors from $i$th population with distribution function $\mathcal{F}_i$, where $\text{E}({\bf x}_{ik}) = {\bs \mu}_i$, $\text{Cov}({\bf x}_{ik}) = {\bs \Sigma}_i > 0$, $i = 1, \ldots, g$, $g \geq 2$, are the mean vector and covariance matrix. Given this set up, we are interested to construct a linear classifier for high-dimensional, low sample size settings, i.e. $p \gg n_i$, and when $\mathcal{F}_i$ can be non-normal. \\
%%%%%%%%%%%%%%%%%%%%%%%%%%%%%%%%%%%%%%%%%%%%%%%%%%%%%%%%%%%%%%%
\indent Classification and regression are two of the most powerful tools of statistical analysis, both as main objective of analysis on their own and also as a source of further investigation. Due to ever growing complexity of data, classification has particularly attracted a central place in modern statistical analysis. The wave of large-dimensional data sets in the last few decades and their associated questions of analysis have lead the researchers to substantially think and improve the classical framework of classification and discrimination. \\
%%%%%%%%%%%%%%%%%%%%%%%%%%%%%%%%%%%%%%%%%%%%%%%%%%%%%%%%%%%%%%%
\indent This paper mainly addresses the classification problem for such a complex data set up, particularly when the dimension of the multivariate vector, $p$, may exceed the number of such vectors, $n_i$, i.e., $p \gg n_i$ (see Sec. \ref{sec:Applns} for examples). As the classical theory of classification does not work in this case, mainly due to the singularity of empirical covariance matrix (see Sec. \ref{sec:TwoSampDF} for more details), efforts have been made in the literature to offer some potential alternatives. \citet{BickelLevina04} propose Independence Rule (IR), or naive Bayes rule, by using only the diagonal of the empirical covariance matrix and compare it to Fisher's linear discriminant function (LDF) for the case of two normal populations. Under certain general conditions on the eigenvalues of the scaled covariance matrix, they show that IR under independence assumption is comparable to Fisher's LDF under dependence when the empirical covariance matrix is replaced with a $g$-inverse computed from the empirical non-zero eigenvalues and the corresponding eigenvectors. A regularized discriminant analysis using Fisher's LDF is given in \citet{WittenTibsh11}. For a short but useful review of this and other classification methods for high-dimensional data, see \citet{Mai13}. Further relevant approaches will be referred to in their appropriate context in the sequel.\\
%%%%%%%%%%%%%%%%%%%%%%%%%%%%%%%%%%%%%%%%%%%%%%%%%%%%%%%%%%%%%%%
\indent We begin in Sec. \ref{sec:TwoSampDF} with the two-sample $U$-classifier, giving detailed explanations on the its construction and justification. An extension to the multi-sample case is given in Sec. \ref{sec:MultSampDF}. Accuracy of the classifier under different parameter settings is shown in Sec. \ref{sec:Simns}, whereas the practical applications on real data sets are demonstrated in Sec. \ref{sec:Applns}. All technical proofs are deferred to the appendix.

%%%%%%%%%%%%%%%%%%%%%%%%%%%%%%%%%%%%%%%%%%%%%%%%%%%%%%%%%%%%%%%%%%%%%%%%%%%%%%%%%%%%%%%%%%%%%%%
\section{The two-sample case}\label{sec:TwoSampDF}
%%%%%%%%%%%%%%%%%%%%%%%%%%%%%%%%%%%%%%%%%%%%%%%%%%%%%%%%%%%%%%%%%%%%%%%%%%%%%%%%%%%%%%%%%%%%%%%

This section is devoted to the construction and justification of the two-sample classifier along with its properties, asymptotic distribution and misclassification rate. These are, respectively, the subjects of the next three subsections.
%%%%%%%%%%%%%%%%%%%%%%%%%%%%%%%%%%%%%%%%%%%%%%%%%%%%%%%%%%%%%%%

%--------------------------------------------------------------------------
\subsection{Construction and motivation of the $U$-classifier}\label{subsec:ConstructDF}
%--------------------------------------------------------------------------

Let ${\bf x}_{ik} = (x_{ik1}, \ldots, x_{ikp})' \sim \mathcal{F}_i$ be as defined above and $\pi_i$ denote the $i$th (unknown) population, $i = 1, 2$ ($g = 2$). Let ${\bf x}_0$ be the new observation to be classified to either of the two populations where the misclassification errors are represented by the conditional probabilities
\[
\pi(i|j) = \text{P}({\bf x}_0 \in \pi_i | \pi_j),~ i, j = 1, 2, ~i \neq j.
\]
Using the information on $p$ characteristics in each sample, i.e. $(x_{ik1}, \ldots, x_{ikp})'$, we aim to construct a classifier which assigns ${\bf x}_0$ to $\pi_i$, $i = 1, 2$, optimally, i.e. by keeping $\pi(i|j)$ as small as possible, as $n_i, p \rightarrow \infty$, when (i) $p$ may arbitrarily exceed $n_i$, $p \gg n_i$, (ii) $\mathcal{F}_i$ may not necessarily be normal, and (iii) ${\bs \Sigma}_i$'s may be unequal, ${\bs \Sigma}_1 \neq {\bs \Sigma}_2$. Note that high-dimensional or, as is frequently known, $(n, p)$-asymptotic framework is kept general in that it implies both $n_i \rightarrow \infty$ and $p \rightarrow \infty$ but without requiring the two indices to satisfy any specific relationship of mutual growth order. It will, however, be shown in the sequel that some of the results hold even by assuming $n_i$ fixed and any arbitrary $p$. \\
%%%%%%%%%%%%%%%%%%%%%%%%%%%%%%%%%%%%%%%%%%%%%%%%%%%%%%%%%%%%%%%
\indent First, to set the notations, let
\begin{equation}\label{eqn:SampMeanCovMat}
\overline{\bf x}_i = \frac{1}{n_i}\sum_{k = 1}^{n_i}{\bf x}_{ik} \quad \text{and} \quad \widehat{\bf \Sigma}_i = \frac{1}{n_i - 1}\sum_{k = 1}^{n_i}({\bf x}_{ik} - \overline{\bf x}_i)({\bf x}_{ik} - \overline{\bf x}_i)
\end{equation}
be the usual unbiased estimators of ${\bs \mu}_i$ and ${\bs \Sigma}_i$, respectively. The classical two-sample linear classifier, assuming equal and known ${\bs \Sigma}$ with equal costs and priors can be expressed, ignoring the constants, as \citep[Ch. 6]{Seber04}
\begin{equation}\label{eqn:ClassTwoSampDF}
A({\bf x}) = (\overline{\bf x}_1 - \overline{\bf x}_2)'{\bs \Sigma}^{-1}{\bf x}_0 - \frac{1}{2}(\overline{\bf x}'_1{\bs \Sigma}^{-1}\overline{\bf x}_1 - \overline{\bf x}'_2{\bs \Sigma}^{-1}\overline{\bf x}_2),
\end{equation}
where ${\bf x}_0$ is the point to be classified. Although this classifier is usually constructed under normality assumption, using a ratio of multivariate normal density functions of the two populations and substituting ${\bs \Sigma}_1 = {\bs \Sigma}_2$ (which makes the classifier linear), it is well-known that Fisher constructed the same classifier without assuming normality, and hence it is also known as Fisher's linear discriminant function. It is the most frequently used classifier in practice and assuming $n_i > p$ and normality, the misclassification probability can be computed using the normal distribution function. Obviously, with ${\bs \Sigma}$ unknown in practice, we need to estimate $A({\bf x})$, replacing ${\bs \Sigma}$ with its usual pooled estimator, $\widehat{\bs \Sigma}_{\text{pooled}} = \sum_{i = 1}^2(n_i - 1)\widehat{\bs \Sigma}_i/\sum_{i = 1}^2(n_i - 1)$, under the homoscedasticity assumption, where $\widehat{\bs \Sigma}_i$ are defined above, so that
\begin{equation}\label{eqn:ClassTwoSampDFEstd}
\widehat{A}({\bf x}) = (\overline{\bf x}_1 - \overline{\bf x}_2)'\widehat{\bs \Sigma}^{-1}_{\text{pooled}}{\bf x}_0 - \frac{1}{2}(\overline{\bf x}'_1\widehat{\bs \Sigma}^{-1}_{\text{pooled}}\overline{\bf x}_1 - \overline{\bf x}'_2\widehat{\bs \Sigma}^{-1}_{\text{pooled}}\overline{\bf x}_2).
\end{equation}
When the data are high-dimensional, i.e., when $p > n_i$, $\widehat{\bs \Sigma}_i$, hence, $\widehat{\bs \Sigma}_{\text{pooled}}$, are singular and can not be inverted, implying that $\widehat{A}({\bf x})$ can not be used in this case. To see how the situation develops in this framework, let us first take $\widehat{\bs \Sigma}$ out of the classifier in (\ref{eqn:ClassTwoSampDF}) and consider
\begin{equation}\label{eqn:ClassTwoSampDFwoSigma}
A({\bf x}) = (\overline{\bf x}_1 - \overline{\bf x}_2)'{\bf x}_0 - \frac{1}{2}(\overline{\bf x}'_1\overline{\bf x}_1 - \overline{\bf x}'_2\overline{\bf x}_2).
\end{equation}
Assuming ${\bf x}_0 \in \pi_1$, we immediately note that
\begin{eqnarray}
\text{E}[A({\bf x})|{\bf x}_0 \in \pi_1]&=& \frac{1}{2}\|{\bs \mu}_1 - {\bs \mu}_2\|^2 - \text{B},\label{eqn:MeanClassDFwBias}
\end{eqnarray}
where $\|\cdot\|$ is the Euclidean norm and B = $\text{tr}({\bs \Sigma}_1)/2n_1 - \text{tr}({\bs \Sigma}_2)/2n_2$. We note that, removing the covariance matrix makes the resulting classifier biased with the bias term, B, composed of the traces of the unknown covariance matrices. If ${\bs \Sigma}_1 = {\bs \Sigma}_2$, then B = $(n_2 - n_1)\text{tr}({\bs \Sigma})/2n_1n_2$, so that the classifier is positively or negatively biased given $n_2 > n_1$ or $n_2 < n_1$. \\
%%%%%%%%%%%%%%%%%%%%%%%%%%%%%%%%%%%%%%%%%%%%%%%%%%%%%%%%%%%%%%%
\indent To inherently adjust the classifier for its bias and improve its accuracy, consider the second component of $A({\bf x})$ in Eqn. (\ref{eqn:ClassTwoSampDFwoSigma}), and note that
\[
\text{E}[(\overline{\bf x}'_1\overline{\bf x}_1 - \overline{\bf x}'_2\overline{\bf x}_2)/2] = \text{B} + ({\bs \mu}'_1{\bs \mu}_1 - {\bs \mu}'_2{\bs \mu}_2)/2,
\]
where ${\bs \mu}'_i{\bs \mu}_i$ were used to complete the squared norm in the expectation of complete classifier in (\ref{eqn:MeanClassDFwBias}) and B is the same bias term. Now
\[
\overline{\bf x}'_i\overline{\bf x}_i  = \frac{1}{n^2_i}\sum_{k = 1}^{n_i}\sum_{r = 1}^{n_i}A_{ikr} = \frac{1}{n^2_i}\sum_{k = 1}^{n_i}A_{ik} + \frac{1}{n^2_i}\underset{k \neq r}{\sum_{k = 1}^{n_i}\sum_{r = 1}^{n_i}}A_{ikr} = Q_{0i} + Q_{1i},
\]
where $A_{ik} = {\bf x}'_{ik}{\bf x}_{ik}$, $A_{ikr} = {\bf x}'_{ik}{\bf x}_{ir}$, $k \neq r$, so that $\text{E}(Q_{0i}) = (\text{tr}({\bs \Sigma}_i) + {\bs \mu}'_i{\bs \mu}_i)/n_i$ and $\text{E}(Q_{1i}) = (1 - 1/n_i){\bs \mu}'_i{\bs \mu}_i$. Let $Q_0 = Q_{01} - Q_{02}$, $Q_1 = Q_{11} - Q_{12}$. Then
\[
\text{E}(Q_0) = 2\text{B} + \text{R} ~~\text{and}~~\text{E}(Q_1) = ({\bs \mu}'_1{\bs \mu}_1 - {\bs \mu}'_2{\bs \mu}_2) - \text{R},
\]
where $R = {\bs \mu}'_1{\bs \mu}_1/n_1 - {\bs \mu}'_2{\bs \mu}_2/n_2$. As R appears with opposite signs in the two components, adjusting each component by this amount, keeping the total expectation same, we have $\text{E}(Q_0)$ + R = ${\bs \mu}'_1{\bs \mu}_1 - {\bs \mu}'_2{\bs \mu}_2$ and $\text{E}(Q_1)$ - R = 2B. Adjusting the corresponding terms in $A({\bf x})$ in (\ref{eqn:ClassTwoSampDFwoSigma}) similarly, we re-write it as
\begin{equation}\label{eqn:ModifiedTwoSDF}
A_0({\bf x}) = \frac{1}{p}(\overline{\bf x}_1 - \overline{\bf x}_2)'{\bf x}_0 - \frac{1}{2}(U_{n_1} - U_{n_2}),
\end{equation}
where $U_{n_i} = \sum_{k \neq r}^{n_i}A_{ikr}/pQ(n_i)$, $Q(n_i) = n_i(n_i - 1)$, is a one-sample $U$-statistic with symmetric kernel, $A_{ikr}/p = {\bf x}'_{ik}{\bf x}_{ir}/p$, $k \neq r$, which is a bilinear form of two independent components. Now, assuming ${\bf x}_0 \in \pi_1$ and independent of the elements of both samples, we have $\text{E}[(\overline{\bf x}_1 - \overline{\bf x}_2)'{\bf x}_0] = {\bs \mu}'_1{\bs \mu}_1 - {\bs \mu}'_1{\bs \mu}_2$, so that
\[
\text{E}[A_0({\bf x})|\pi_1] = \|{\bs \mu}_1 - {\bs \mu}_2\|^2/2,
\]
without any bias term. Further, with ${\bf x}_0 \in \pi_1$, $A_0({\bf x})$ is composed of four bilinear forms, two from sample 1, one from sample 2, and one mixed. By symmetry, $\text{E}[A_0({\bf x})|\pi_2] = -\|{\bs \mu}_1 - {\bs \mu}_2\|^2/2$, and the classifier again consists of four bilinear forms, two from sample 2, one from sample 1 and one mixed. We, therefore, define the classification rule for the proposed $U$-classifier as
%%%%%%%%%%%%%%%%%%%%%%%%%%%%%%%%%%%%%%%%%%%%%%%%%%%%%%%%%%%%%%%
\[
\text{Assign ${\bf x}_0$ to $\pi_1$ if $A_0({\bf x}) > 0$, otherwise to $\pi_2$.}
\]
%%%%%%%%%%%%%%%%%%%%%%%%%%%%%%%%%%%%%%%%%%%%%%%%%%%%%%%%%%%%%%%
Before we study the properties of $A_0({\bf x})$ in the next section, a few important remarks are in order. First, $A_0({\bf x})$ is composed of bilinear forms - and we call it \emph{bilinear classifier} - where the bi-linearity of the $U$-component is expressed in the kernels of the two $U$-statistics and that of the $P$-component by the projection of the new observation with respect to the difference between the empirical centroids of the two independent samples. Further, $A_0({\bf x})$ is entirely composed of empirical quantities, free of any unknown parameter, so that it can be directly used in practice using the decision rule stated above. Note also that, $A_0({\bf x})$ is linear but the linearity does not require homoscedasticity assumption, i.e., it is linear even if ${\bs \Sigma}_1 \neq {\bs \Sigma}_2$.\\
%%%%%%%%%%%%%%%%%%%%%%%%%%%%%%%%%%%%%%%%%%%%%%%%%%%%%%%%%%%%%%%
\indent Moreover, the first part of $A_0({\bf x})$ is normalized by $p$, and so are the kernels of $U$-statistics in the second part. This will help us derive the limiting distribution of the classifier for $(n, p)$-asymptotics under a general multivariate model and mild assumptions. As a final remark, recall that the formulation of $A_0({\bf x})$ arises from depriving the original classifier of empirical covariance matrix. Whereas, this removal of an essential ingredient has its price to be paid, the resulting classifier can still be justified from a different perspective which has its merit, particularly for high-dimensional data.\\
%%%%%%%%%%%%%%%%%%%%%%%%%%%%%%%%%%%%%%%%%%%%%%%%%%%%%%%%%%%%%%%
\indent To see this, consider $\overline{d}_{12} = \overline{d}_1 - \overline{d}_2$, $\overline{d}_i = \sum_{k = 1}^{n_i}d_{ki}/n_i$, where $d_{ki} = \|{\bf x}_0 - {\bf x}_{ki}\|^2$ is the Euclidean distance of ${\bf x}_0$ from sample $i = 1, 2$, $k = 1, \ldots, n_i$. It follows that $\overline{d}_{12}$ has same bias B as for $A({\bf x})$. For all expressions location-invariant, write $\widehat{\bs \Sigma}_i = \sum_{k \neq r}^{n_i}{\bf D}'_{ikr}{\bf D}_{ikr}/n_i(n_i - 1)$, ${\bf D}_{ikr} = {\bf x}_{ik} - {\bf x}_{ir}$ \citep{Ahmad14b}. Now $\widetilde{A}({\bf x}) = \overline{d}_{12} - [\text{tr}(\widehat{\bs \Sigma}_1)/n_1 - \text{tr}(\widehat{\bs \Sigma}_2)/n_2]$, and since $\sum_{k = 1}^{n_i}d_i = (n_i - 1)\text{tr}(\widehat{\bs \Sigma}_i) + n_i\|{\bf x}_0 - \overline{\bf x}_i\|^2$, it simplifies to $\widetilde{A}({\bf x}) = A({\bf x})$ + B; compare with Eqn. (\ref{eqn:MeanClassDFwBias}).\\
%%%%%%%%%%%%%%%%%%%%%%%%%%%%%%%%%%%%%%%%%%%%%%%%%%%%%%%%%%%%%%%
\indent This implies that $A_0({\bf x})$ can also be constructed using Euclidean distances and the same bias-adjustment that lead Eqn. (\ref{eqn:MeanClassDFwBias}) to Eqn. (\ref{eqn:ModifiedTwoSDF}). This distance-based approach has been discussed in \citet{ChanHall09} and the same is further evaluated in \citet{AoshimaYata14}. Our approach, however, makes the classifier not only unbiased but also more general as well as convenient to study and apply in practice.

%--------------------------------------------------------------------------
\subsection{Asymptotic distribution of the $U$-classifier}\label{subsec:AsympDistnDF}
%--------------------------------------------------------------------------

Given ${\bf x}_{ik} \sim \mathcal{F}_i$, let ${\bf z}_{ik} = {\bf x}_{ik} - {\bs \mu}_i$ with $\text{E}({\bf z}_{ik}) = {\bf 0}$, $\text{Cov}({\bf z}_{ik}) = {\bs \Sigma}_i$, $i = 1, 2$. When we relax normality, we assume the following multivariate model
\begin{equation}\label{eqn:ModelLINN}
{\bf z}_{ik} = {\bs \Lambda}_i{\bf y}_{ik},
\end{equation}
where ${\bf y}_{ik} = (y_{ik1}, \ldots, y_{ikp})'$ has iid elements with $\text{E}({\bf y}_{ik}) = {\bf 0}$, $\text{Cov}({\bf y}_{ik}) = {\bf I}$, and ${\bs \Lambda}_i$ is a known $p\times p$ matrix of constants such that ${\bs \Lambda}'_i{\bs \Lambda}_i = {\bf A}_i$, ${\bs \Lambda}_i{\bs \Lambda}'_i = {\bs \Sigma}_i > 0$, $i = 1, 2$. To study the properties of the classifier and its asymptotic normality, we shall supplement Model (\ref{eqn:ModelLINN}) with the following assumptions.
%------------------------------------------------------------%
\begin{assn}\label{assn:4thMomnt}
$\text{E}(y^4_{iks}) = \gamma < \infty$, $\gamma \in \mathbb{R}^+$, $i = 1, 2$.
\end{assn}
\begin{assn} \label{assn:traceSigma}
$\lim_{p \rightarrow \infty}\text{tr}({\bs \Sigma}_i)/p = O(1)$, $i = 1, 2$.
\end{assn}
\begin{assn}\label{assn:ExtraH1Distn}
$\lim_{p \rightarrow \infty}{\bs \mu}'_i{\bs \Sigma}_k{\bs \mu}_j/p = O(1)$, $i, j , k = 1, 2$.%, $k = i$ or $k = j$, $i = j$  or $i \neq j$.
\end{assn}
\begin{assn}\label{assn:traceSigmaHadProd}
$\lim_{p \rightarrow \infty}\frac{\text{tr}({\bs \Sigma}^a_i\odot{\bs \Sigma}^b_j)}{\text{tr}({\bs \Sigma}^a_i\otimes{\bs \Sigma}^b_j)} = 0$, $a, b = 1, 2, 3$, $a + b \leq 4$, $i, j = 1, 2$, where $\odot$ and $\otimes$ are Hadamard and Kronecker products.
\end{assn}
%--------------------------------------------------------------%
Assumption \ref{assn:4thMomnt} essentially replaces normality. Assumptions \ref{assn:traceSigma} is simple and mild, and as its consequence, $\lim_{p \rightarrow \infty}\text{tr}({\bs \Sigma}_i{\bs \Sigma}_j)/p^2 = O(1)$, so that a reference to the assumption may also imply it consequence in the sequel. Assumptions \ref{assn:ExtraH1Distn} and \ref{assn:traceSigmaHadProd} are needed only to show control of misclassification rate and consistency of the moments of classifier. Assumption \ref{assn:traceSigmaHadProd} ensures that the moments asymptotically coincide with those under normality. This assumption is neither directly needed in practical use of the classifier, nor is it required under normality whence all terms involving the ratio vanish. The same assumptions will be extended for multi-sample case in Sec. \ref{sec:MultSampDF}. \\
%%%%%%%%%%%%%%%%%%%%%%%%%%%%%%%%%%%%%%%%%%%%%%%%%%%%%%%%%%%%%%%
\indent Now, continuing to assume ${\bf x}_0$ independent of all elements of sample 1 (where it is already independent of all elements of sample 2), we get the following lemma, proved in Appendix \ref{subsec:ProofLemmaMomsTwoSDF}, on the moments of classifier.
%%%%%%%%%%%%%%%%%%%%%%%%%%%%%%%%%%%%%%%%%%%%%%%%%%%%%%%%%%%%%%%
%%%%%%%%%%%%%%%%%%%%%%%%%%%%%%%%%%%%%%%%%%%%%%%%%%%%%%%%%%%%%%%
\begin{lem}\label{lem:MomsTwoSDF} Let the two-sample modified classifier be as given in Eqn. (\ref{eqn:ModifiedTwoSDF}). Then, assuming ${\bf x}_0 \in \pi_i$, we have
\begin{eqnarray}
\text{E}[A_0({\bf x})|\pi_i] &=& \frac{1}{2p}\|{\bs \mu}_i - {\bs \mu}_j\|^2_{\bf I} = \frac{1}{2p}\Delta^2_{\bf I}\label{eqn:MeanDFTwoS}\\
\text{Var}[A_0({\bf x})|\pi_i] &=& \frac{\delta^2_i}{p^2} + \frac{1}{p^2}\|{\bs \mu}_1 - {\bs \mu}_2\|^2_{{\bs \Sigma}^{-1}_i} = \frac{\delta^2_i}{p^2} + \frac{1}{p^2}\Delta^2_{{\bs \Sigma}^{-1}_i},\label{eqn:VarDFTwoS}
\end{eqnarray}
where $\Delta^2_{{\bf M}} = \|{\bs \mu}_i - {\bs \mu}_j\|^2_{\bf M} = ({\bs \mu}_i - {\bs \mu}_j)'{\bf M}^{-1}({\bs \mu}_i - {\bs \mu}_j)$, $i, j = 1, 2$, $i \neq j$, for any ${\bf M}_p > 0$, and
\[
\delta^2_i = \frac{\text{tr}({\bs \Sigma}^2_i)}{n_i} + \frac{\text{tr}({\bs \Sigma}_i{\bs \Sigma}_j)}{n_j} + \sum_{i = 1}^2\frac{\text{tr}({\bs \Sigma}^2_i)}{2n_i(n_i - 1)}.
\]
\end{lem}
%%%%%%%%%%%%%%%%%%%%%%%%%%%%%%%%%%%%%%%%%%%%%%%%%%%%%%%%%%%%%%%
%%%%%%%%%%%%%%%%%%%%%%%%%%%%%%%%%%%%%%%%%%%%%%%%%%%%%%%%%%%%%%%
The moments in Lemma \ref{lem:MomsTwoSDF} are reported using general notation for ${\bf x}_0 \in \pi_i$ so that they can be easily extended to the multi-sample case later. For the present case of $g = 2$, the moments easily reduce to their specific form for $i = 1$ or $i = 2$, where the mean is obviously the same for both, i.e., $\|{\bs \mu}_1 - {\bs \mu}_1\|^2_{\bf I}/2p$.\\
%%%%%%%%%%%%%%%%%%%%%%%%%%%%%%%%%%%%%%%%%%%%%%%%%%%%%%%%%%%%%%%
\indent With the picture of the proposed classifier and its moments relatively clear, we can express our high-dimensional classification problem precisely as
\begin{equation}\label{eqn:ClassfcnProb}
\mathcal{C}_p = \{({\bs \mu}_i, {\bs \Sigma}_i):~\Delta^2_{\bf I}, \Delta^2_{{\bs \Sigma}^{-1}_i}\}_p,~i = 1, 2,
\end{equation}
where the index $p$ implies the dependence of components $\mathcal{C}_p$ on the dimension. Note that, the second component in Eqn. (\ref{eqn:VarDFTwoS}) vanishes under Assumption \ref{assn:ExtraH1Distn}. The rest of $\text{Var}[A_0({\bf x})|\pi_i]$, and $\text{E}[A_0({\bf x})|\pi_i]$, are uniformly bounded in $p$, for any fixed $n_i$, under Assumptions \ref{assn:traceSigma}. We can thus write
\begin{eqnarray}
\lim_{p \rightarrow \infty}\text{E}[A_0({\bf x})|\pi_i] &=& \frac{1}{2}\Delta^2_{0, {\bf I}} \label{eqn:LmtgMeanDF}\\
\lim_{p \rightarrow \infty}\text{Var}[A_0({\bf x})|\pi_i] &=& \delta^2_{0, i}\left[O\left(\frac{1}{n_1} + \frac{1}{n_2}\right) + o(1)\right],\label{eqn:LmtgVarDF}
\end{eqnarray}
where $\Delta^2_{0, {\bf I}} = \lim_{p \rightarrow \infty} \Delta^2_{\bf I}/p \in (0, \infty)$ and $\delta^2_{0, _i} = \lim_{p \rightarrow \infty}\delta^2_i/p^2 \in (0, \infty)$. Now, the variance obviously vanishes when we also allow $n_i \rightarrow \infty$ along with $p \rightarrow \infty$, which immediately gives consistency of the classifier, formally stated in Theorem \ref{thm:ConsistencyTwoSDF} below and proved in Appendix \ref{subsec:ProofThmConsisTWoSDF}. Obviously, in practice, the consistency is expected to hold with unknown parameters replaced by their estimators. We need to estimate $\Delta^2_{\bf I}$ and non-vanishing traces in $\delta^2_i$ to estimate the limiting moments of the classifier. In the following, we define unbiased and consistent plug-in estimators of these components.\\
%%%%%%%%%%%%%%%%%%%%%%%%%%%%%%%%%%%%%%%%%%%%%%%%%%%%%%%%%%%%%%%
\indent For ${\bf x}_{ik}$, $i = 1, 2$, let ${\bf D}_{ikr} = {\bf x}_{ik} - {\bf x}_{ir}$, $k \neq r$, be as defined in Sec. \ref{subsec:ConstructDF} with $\text{E}({\bf D}_{ikr}) = 0$, $\text{Cov}({\bf D}_{ikr}) = 2{\bs \Sigma}_i$. Also let ${\bf D}_{ijkl} = {\bf x}_{ik} - {\bf x}_{jl}$ using differences of vectors from two independent samples with $\text{E}({\bf D}_{ijkl}) = {\bs \mu}_i - {\bs \mu}_j$, $\text{Cov}({\bf D}_{ijkl}) = {\bs \Sigma}_i + {\bs \Sigma}_j$, $i, j = 1, 2$, $i \neq j$. Extending the strategy of Sec. \ref{subsec:ConstructDF} for the estimation of $\text{tr}({\bs \Sigma}_i)$, we define estimators of the traces involved in $\delta^2_i$. Let $A^2_{ikrk'r'} = ({\bf D}'_{ikr}{\bf D}_{ik'r'})^2$ and $A^2_{ijkrls} = ({\bf D}'_{ikr}{\bf D}_{jls})^2$, $A_{ijkl} = {\bf D}'_{ik}{\bf D}_{jl}$, $i \neq j$, using within- and between-sample independent vectors, respectively. Then, by independence, the plug-in estimators of $\Delta^2_{\bf I}/p^2$, $\text{tr}({\bs \Sigma}^2_i)/p^2$ and $\text{tr}({\bs \Sigma}_i{\bs \Sigma}_j)/p^2$ are defined, respectively, as
\begin{eqnarray}
E_0 &=& U_{n_i} + U_{n_j} - 2U_{n_in_j}\label{eqn:LocInvE0}\\
E_i &=& \frac{1}{12\eta(n_i)}\underset{\pi(k, r, k', r')}{\sum_{k = 1}^{n_i}\sum_{r = 1}^{n_i}\sum_{k' = 1}^{n_i}\sum_{r' = 1}^{n_i}}\frac{1}{p^2}A^2_{ikrk'r'}\label{eqn:LocInvE1}\\
E_{ij} &=& \frac{1}{4\eta(n_in_j)}\underset{\pi(k, r)}{\sum_{k = 1}^{n_i}\sum_{r = 1}^{n_i}}\underset{\pi(l, s)}{\sum_{l = 1}^{n_j}\sum_{s = 1}^{n_j}}\frac{1}{p^2}A^2_{ijklrs}\label{eqn:LocInvE2}
\end{eqnarray}
where $\eta(n_i) = n_i(n_i - 1)(n_i - 2)(n_i - 3)$ and $\eta(n_in_j) = n_in_j(n_i - 1)(n_j - 1)$, $i, j = 1, 2$, $i \neq j$, and $\pi(\cdot)$ implies all indices pairwise unequal. Further,
\[
U_{n_in_j} = \frac{1}{n_in_j}\sum_{k = 1}^{n_i}\sum_{l = 1}^{n_j}\frac{1}{p}A_{ijkl},
\]
a two-sample $U$-statistic, so that, with $U_{n_i}$, $i = 1, 2$, as one-sample $U$-statistic defined after Eqn. (\ref{eqn:ModifiedTwoSDF}), $E_0$ estimates $\Delta^2_{\bf I}/p^2 = \|{\bs \mu}_i - {\bs \mu}_j\|^2/p^2$. $E_0$, $E_i$, $E_{ij}$ are unbiased and location-invariant estimators. The following theorem, proved in Appendix \ref{subsec:ProofThmPropsEstrs}, shows further that the variances of the ratios of these estimators to the parameters they estimate are uniformly bounded in $p$, so that they are consistent as $p \rightarrow$ ($n_i$ fixed) and also when $n_i, p \rightarrow \infty$.
%%%%%%%%%%%%%%%%%%%%%%%%%%%%%%%%%%%%%%%%%%%%%%%%%%%%%%%%%%%%%
%%%%%%%%%%%%%%%%%%%%%%%%%%%%%%%%%%%%%%%%%%%%%%%%%%%%%%%%%%%%%
\begin{thm}\label{thm:PropsEstrsBounds} $E_0, E_i$, $E_{ij}$, defined in Eqns. (\ref{eqn:LocInvE0})-(\ref{eqn:LocInvE2}), are unbiased estimators of $\Delta^2_{\bf I}/p^2$, $\text{tr}({\bs \Sigma}^2_i)/p^2$ and $\text{tr}({\bs \Sigma}_1{\bs \Sigma}_2)/p^2$. Further, under Assumptions \ref{assn:4thMomnt}-\ref{assn:traceSigmaHadProd},
\begin{eqnarray}
\text{Var}\left(\frac{E_0}{\Delta^2_{\bf I}}\right) &=& O\left(\frac{1}{n_i} + \frac{1}{n_j}\right)\label{eqn:BoundVarE0}\\
\text{Var}\left(\frac{E_i}{\text{tr}({\bs \Sigma}^2_i)}\right) &=& O\left(\frac{1}{n_i}\right)\label{eqn:BoundVarE1}%\\
\end{eqnarray}
\begin{eqnarray}
\text{Var}\left(\frac{E_{ij}}{\text{tr}({\bs \Sigma}_i{\bs \Sigma}_j)}\right) &=& O\left(\frac{1}{n_i} + \frac{1}{n_j}\right)\label{eqn:BoundVarE2}\\
\text{Cov}\left(\frac{E_i}{\text{tr}({\bs \Sigma}^2_i)}, \frac{E_{ij}}{\text{tr}({\bs \Sigma}_i{\bs \Sigma}_j)}\right) &=& O\left(\frac{1}{n_i}\right).\label{eqn:BoundCovs}
\end{eqnarray}
\end{thm}
%%%%%%%%%%%%%%%%%%%%%%%%%%%%%%%%%%%%%%%%%%%%%%%%%%%%%%%%%%%%%
%%%%%%%%%%%%%%%%%%%%%%%%%%%%%%%%%%%%%%%%%%%%%%%%%%%%%%%%%%%%%
The bounds in Theorem \ref{thm:PropsEstrsBounds} suffice for consistency of estimators so that exact variances and covariances are not needed. These exact moments follow from Theorem \ref{thm:BasicQFBFMomsNN} and Lemma \ref{lem:QFBFResults}; see \citet{Ahmad14b}.\\
%%%%%%%%%%%%%%%%%%%%%%%%%%%%%%%%%%%%%%%%%%%%%%%%%%%%%%%%%%%%%
\indent From Theorem \ref{thm:PropsEstrsBounds}, it immediately follow that $E_0/\text{E}(E_0) \xrightarrow{\mathcal{P}} 1$ which gives empirical mean, $\widehat{\text{E}}[A_0({\bf x})|\pi_i] = \frac{1}{2}E_0$, as a consistent estimator of the true mean of the classifier. Using similar probability convergence of the other two estimators $E_i$ and $E_{ij}$, we obtain, by Slutsky's lemma \citep[][p 11]{vdv98}, the first component of the empirical variance, $\widehat{\text{Var}}[A_0({\bf x})|\pi_i]$, i.e. $\widehat{\delta}^2_i$, as a consistent estimator of $\delta^2_i$, $i = 1, 2$, using plug-in estimators such that $\widehat{\delta}^2_i/\delta^2_i \xrightarrow{\mathcal{P}} 1$. The limiting empirical moments, parallel to Eqns. (\ref{eqn:LmtgMeanDF})-(\ref{eqn:LmtgVarDF}), thus follow as
\begin{eqnarray}
\lim_{p \rightarrow \infty}\widehat{\text{E}}[A_0({\bf x})|\pi_i] &=& \frac{1}{2}\Delta^2_{0, {\bf I}}\label{eqn:LmtgMeanDFEstd}\\
\lim_{p \rightarrow \infty}\widehat{\text{Var}}[A_0({\bf x})|\pi_i] &=& \delta^2_{0, i}\left[O\left(\frac{1}{n_1} + \frac{1}{n_2}\right) + o_P(1)\right]\label{eqn:LmtgVarDFEstd}
\end{eqnarray}
with $\Delta^2_{0, {\bf I}} = \lim_{p \rightarrow \infty} E_0/p \in (0, \infty)$, $\delta^2_{0, i} = \lim_{p \rightarrow \infty}\widehat{\delta}^2_i/p^2 \in (0, \infty)$. Hence
\begin{eqnarray}
\lim_{p \rightarrow \infty}\left[\widehat{\text{E}}[A_0({\bf x})|\pi_i] - \text{E}[A_0({\bf x})|\pi_i]\right] &=& o_P(1)\label{eqn:LmtgMeanDFDiff}\\
\lim_{p \rightarrow \infty}\left[\widehat{\text{Var}}[A_0({\bf x})|\pi_i] - \text{Var}[A_0({\bf x})|\pi_i]\right] &=&  o_P(1).\label{eqn:LmtgVarDFDiff}
\end{eqnarray}
The following theorem, proved in Appendix \ref{subsec:ProofThmConsisTWoSDF}, summarizes both true and empirical consistency of the classifier.
%%%%%%%%%%%%%%%%%%%%%%%%%%%%%%%%%%%%%%%%%%%%%%%%%%%%%%%%%%%%%%%
%%%%%%%%%%%%%%%%%%%%%%%%%%%%%%%%%%%%%%%%%%%%%%%%%%%%%%%%%%%%%%%
\begin{thm}\label{thm:ConsistencyTwoSDF} Given $A_0({\bf x})$ in Eqn. (\ref{eqn:ModifiedTwoSDF}) with its moments as in Lemma \ref{lem:MomsTwoSDF}. Let ${\bf x}_0 \in \pi_i$. Under Assumptions \ref{assn:4thMomnt}-\ref{assn:ExtraH1Distn}, as $n_i, p \rightarrow \infty$, $i = 1, 2$,
\[
\frac{A_0({\bf x})}{\Delta^2_{{\bf I}}/p} \xrightarrow{\mathcal{P}} \frac{(-1)^i}{2} + o_P(1),
\]
with $\delta_i$ defined above. Further, the consistency holds when the moments of the classifier are replaced with their empirical estimators, given in Eqns. (\ref{eqn:LmtgMeanDFEstd})-(\ref{eqn:LmtgVarDFEstd}).
\end{thm}
%%%%%%%%%%%%%%%%%%%%%%%%%%%%%%%%%%%%%%%%%%%%%%%%%%%%%%%%%%%%%%%
%%%%%%%%%%%%%%%%%%%%%%%%%%%%%%%%%%%%%%%%%%%%%%%%%%%%%%%%%%%%%%%
The same arguments help us establish asymptotic normality of the classifier as stated in the following theorem, proved in Appendix \ref{subsec:ProofThmAsympNTwoSDF}.
%%%%%%%%%%%%%%%%%%%%%%%%%%%%%%%%%%%%%%%%%%%%%%%%%%%%%%%%%%%%%%%
%%%%%%%%%%%%%%%%%%%%%%%%%%%%%%%%%%%%%%%%%%%%%%%%%%%%%%%%%%%%%%%
\begin{thm}\label{thm:AsympNTwoSDF} Given $A_0({\bf x})$ in Eqn. (\ref{eqn:ModifiedTwoSDF}) with its moments as in Lemma \ref{lem:MomsTwoSDF}. Let ${\bf x}_0 \in \pi_i$. Under Assumptions \ref{assn:4thMomnt}-\ref{assn:ExtraH1Distn}, as $n_i, p \rightarrow \infty$, $i = 1, 2$,
\[
\frac{A_0({\bf x}) - \text{E}[A_0({\bf x})]}{\sqrt{\text{Var}[A_0({\bf x})]}} \xrightarrow{\mathcal{D}} N(0, 1),
\]
Further, the normal limit holds when the moments are replaced with their empirical estimators, given in Eqns. (\ref{eqn:LmtgMeanDFEstd})-(\ref{eqn:LmtgVarDFEstd}).
\end{thm}
%%%%%%%%%%%%%%%%%%%%%%%%%%%%%%%%%%%%%%%%%%%%%%%%%%%%%%%%%%%%%%%
%%%%%%%%%%%%%%%%%%%%%%%%%%%%%%%%%%%%%%%%%%%%%%%%%%%%%%%%%%%%%%%
The construction of $A_0({\bf x})$ is of great benefit in proving Theorem \ref{thm:AsympNTwoSDF}. Its composition of two parts, each of which in turn a linear combination of two independent components, reduces the bulk of computational burden. Moreover, the optimality property of $U$-statistics ensures the minimum variance (efficiency) of the classifier. A further verification of these properties of the classifier through simulations is demonstrated in Sec. \ref{sec:Simns}.

%%%%%%%%%%%%%%%%%%%%%%%%%%%%%%%%%%%%%%%%%%%%%%%%%%%%%%%%%%%%%%%%%%%%%%%%%%%%%%%%%%%%%%%%%%%%%%%
\section{Estimation of misclassification probabilities}\label{sec:MisclassfnProbs}
%%%%%%%%%%%%%%%%%%%%%%%%%%%%%%%%%%%%%%%%%%%%%%%%%%%%%%%%%%%%%%%%%%%%%%%%%%%%%%%%%%%%%%%%%%%%%%%

Consider the classification problem $\mathcal{C}_p$ in (\ref{eqn:ClassfcnProb}) again. Using notation introduced in the beginning of Sec. \ref{subsec:ConstructDF}, the optimality of the proposed classifier can be evaluated by the misclassification rates $\pi(1|2)$ and $\pi(2|1)$ with
\begin{equation}
\pi(i|j) = P({\bf x} \in \pi_i | \pi_j) = \int_{\mathbb{R}_i}d\mathcal{F}_j({\bf x}),\label{eqn:MisclRate}
\end{equation}
where $\mathcal{F}_j$ denotes the distribution function and $\mathbb{R}_i = \{{\bf x}: {\bf x} \in \pi_i\}$ is the region of observed data from $i$th population with $\mathbb{R}_1 \cup \mathbb{R}_2 = \mathcal{X}$, $\mathbb{R}_1 \cap \mathbb{R}_2 = \emptyset$, where $\mathcal{X}$ denotes the space of observed ${\bf x}$ and $\emptyset$ is the empty set. Under the assumption of equal probabilities and equal costs, we are interested to minimize the total probability of misclassification, say $\epsilon_{o}$, given the observed data, i.e.
\[
\arg\min_{\bf x} \epsilon_0 = \arg\min_{\bf x} [\pi(1|2) + \pi(2|1)]/2,
\]
where the subscript $o$ stands for \emph{optimal}. Obviously, the ideal minimum can only be achieved when the parameters are known in which case the (ideal) classifier takes the form
\[
A^{\text{ideal}}_0({\bf x} \in \pi_1) = {\bf x}'_0({\bs \mu}_1 - {\bs \mu}_2)/p - ({\bs \mu}'_1{\bs \mu}_1 - {\bs \mu}'_2{\bs \mu}_2)/2p,
\]
using the fact that $\overline{\bf x}_i$ and $U_{n_i}$ are unbiased estimators of ${\bs \mu}_i$ and ${\bs \mu}'_i{\bs \mu}_i$, respectively, $i = 1, 2$. Now, if $\mathcal{F}_i$ are known, say multivariate normal, i.e. ${\bf x}_{ik} \sim \mathcal{N}_p({\bs \mu_{i}}, {\bs \Sigma}_{i})$, then, under the homoscedasticity assumption ${\bs \Sigma}_1 = {\bs \Sigma}_2 = {\bs \Sigma}$, the error rate of $A^{\text{ideal}}_0$  can be expressed as
\[
\epsilon_o^{\text{ideal}} = {\bs \Phi}\left(-\frac{ \|{\bs \mu}_1 - {\bs \mu}_2 \|^{2}_{\bf I}}
{2\sqrt{ \|{\bs \mu}_1 - {\bs \mu}_2 \|^{2}_{\bf \Sigma^{-1}} } }  \right),
\]
where ${\bs \Phi}$ is the standard normal distribution function. Assuming equal priors, the best possible performance in this ideal setting, i.e. with ${\bs \mu}_1$, ${\bs \mu}_2$, ${\bf \Sigma}$ known, is achieved by Fisher's linear classifier (equivalently, Bayes rule)
\[
A^{\text{Fisher}}_0 ({\bf x}) = {\bf x}'_0 {\bf \Sigma}^{-1}({\bs \mu}_1 - {\bs \mu}_2) - \frac{1}{2} \left ({\bs \mu}_1 + {\bs \mu}_2 \right )' {\bf \Sigma}^{-1} ({\bs \mu}_1 - {\bs \mu}_2)
\]
\citep{Anderson03} with the corresponding misclassification rate given by
\[
\epsilon^{\textrm{Fisher}}_o = {\bs \Phi} \left(-\frac{1}{2}
\sqrt{\| {\bs \mu}_1 - {\bs \mu}_2 \|^{2}_{\bf \Sigma}} \right),
\]
where $\| {\bs \mu}_1 - {\bs \mu}_2 \|^{2}_{\bf \Sigma}$ is the Mahalanobis distance. Denoting $\epsilon^{\text{Fisher}}_o$ as a benchmark, the relative performance of $A^{\text{ideal}}_0({\bf x})$ can be theoretically evaluated by using the ratio of the arguments of ${\bs \Phi}$, say $q$, where
\[
q = \frac{\|{\bs \mu}_1 - {\bs \mu}_2 \|^{2}_{\bf I} } { \left [ \| {\bs \mu}_1 - {\bs \mu}_2 \|^{2}_{\bf \Sigma} \cdot \| {\bs \mu}_1 - {\bs \mu}_2 \|^{2}_{\bf \Sigma^{-1}}   \right ]^{1/2}}.
\]
\citet{BickelLevina04} put forth a nice strategy to compute a bound for an expression like $q$, based on Kantorovich inequality \citep{Bernstein09}. Following the same idea, let $\bf M$ be any positive definite symmetric $p \times p$ matrix. Then for any vector ${\bs \upsilon}$
\[
\frac{\|{\bs \upsilon} \|^{2}_{\bf I}  } { \|{\bs \upsilon} \|^{2}_{{\bf M}} \cdot \|{\bs \upsilon} \|^{2}_{{\bf M}^{-1}} } \geq \frac{4 \lambda_{\min}(\bf M) \cdot  \lambda_{\max}({\bf M})}
{\left [\lambda_{\min}({\bf M})+\lambda_{\max}({\bf M})\right ]^{2} },
 \]
where $\lambda_{\min}({\bf M})$ and $\lambda_{\max}({\bf M})$ denote the smallest and the largest eigenvalues of  ${\bf M}$, respectively. Applying this inequality to $q$ and denoting the ratio $\frac{\lambda_{\max}({\bf \Sigma})}{\lambda_{\min}({\bf \Sigma})} = \kappa$ (assuming the two extreme eigenvalues bounded away from $0$ and $\infty$), we get
\begin{equation}\label{eqn:BoundQ}
q \geq \frac{2\sqrt{\kappa}}{1+\kappa},
\end{equation}
so that the upper bound on the misclassification probability of $A^{\text{ideal}}_0({\bf x})$ is
\[
\epsilon_o^{\text{ideal}} \leq {\bs \Phi} \left ( - \frac{2\sqrt{\kappa}}{1+\kappa} {\bs \Phi}^{-1} \left(- \epsilon^{\text{Fisher}}_o \right ) \right ).
\]
which essentially depends on $\kappa$, the range of non-zero eigenvalues of ${\bf \Sigma}$.
%%%%--------------------------------------------------------%%%
\begin{figure}[t!]\centering\label{fig:BoundQ}
\vspace{-2cm}	
%\begin{center}
\includegraphics[width = 9.5cm]{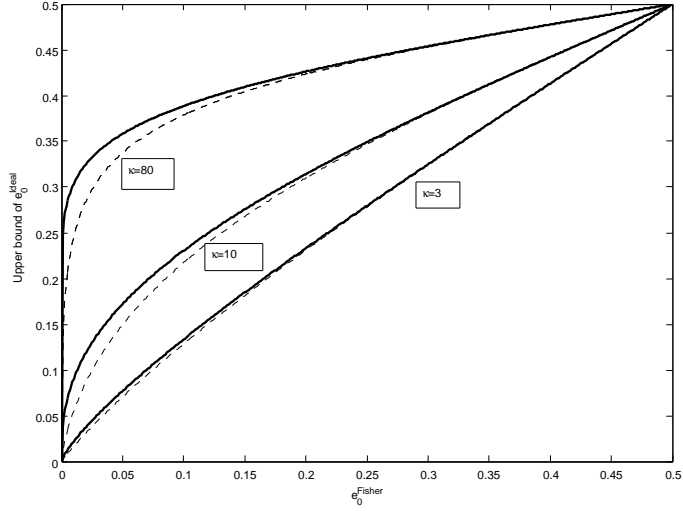}\vspace{-2cm}
\caption{\footnotesize \sc Upper bound on the misclassification probability of $A^{\text{ideal}}_0({\bf x})$ as a function of $\epsilon^{\text{Fisher}}_o$ for normal (thick line) and $t_5$ (dashed line) distributions with $\kappa = 3, 10, 80$.}
%\end{center}
\end{figure}
%%%%--------------------------------------------------------%%%
%%%%%%%%%%%%%%%%%%%%%%%%%%%%%%%%%%%%%%%%%%%%%%%%%%%%%%%%%%%%%%%
We note that, for moderate $\kappa$, the increase in the misclassification rate, induced by taking the  covariance matrix away while constructing $A^{\text{ideal}}_0({\bf x})$, is not large relative to the best possible performance, i.e. $\epsilon_o^{\text{Fisher}}$ (see Fig. 1). Further, the upper bound in (\ref{eqn:BoundQ}) represents the worst-case scenario so that the empirical results are expected to be better.\\
%%%%%%%%%%%%%%%%%%%%%%%%%%%%%%%%%%%%%%%%%%%%%%%%%%%%%%%%%%%%%%%
\indent Now, for an alternative flavor of the evaluation of the classifier, while still continuing to assume normality, let us condition the classifier on the data, i.e., $\underline{A}_0({\bf x}) = A_0({\bf x})|(\overline{\bf x}_i, U_{n_i})$, say, and it immediately follows that
\[
\underline{A}_0({\bf x}) \sim N\left(\frac{1}{p}{\bs \mu}'_i(\overline{\bf x}_1 - \overline{\bf x}_2) - \frac{1}{2}(U_{n_1} - U_{n_2}), ~\|\overline{\bf x}_1 - \overline{\bf x}_2\|^2_{{\bs \Sigma}^{-1}_i}\right)
\]
$i = 1, 2$. This, using the standardized version of classifier (Theorem \ref{thm:AsympNTwoSDF}), leads to the actual error rate
\begin{eqnarray*}
\underline{\epsilon}_n &=& \frac{1}{2}\left[\Phi\left(-\frac{{\bs \mu}'_1(\overline{\bf x}_1 - \overline{\bf x}_2) - (U_{n_1} - U_{n_2})/2}{\sqrt{(\overline{\bf x}_1 - \overline{\bf x}_2)'{\bs \Sigma}_1(\overline{\bf x}_1 - \overline{\bf x}_2)}}\right)\right.\\
&&\qquad\qquad\quad\qquad\quad+~\left.\Phi\left(-\frac{{\bs \mu}'_2(\overline{\bf x}_2 - \overline{\bf x}_1) - (U_{n_1} - U_{n_2})/2}{\sqrt{(\overline{\bf x}_1 - \overline{\bf x}_2)'{\bs \Sigma}_2(\overline{\bf x}_1 - \overline{\bf x}_2)}}\right)\right],
\end{eqnarray*}
where the subscript $n$ denotes the dependence on the observed sample. Using Theorem \ref{thm:ConsistencyTwoSDF},
\begin{eqnarray*}
\left|\left\{\frac{1}{p}{\bs \mu}'_i(\overline{\bf x}_1 - \overline{\bf x}_2) - \frac{1}{2}(U_{n_1} - U_{n_2})\right\} - \frac{1}{2}\Delta^2_{0, {\bf I}}\right| &\xrightarrow{\mathcal{P}}& 0\\
\left|\|\overline{\bf x}_1 - \overline{\bf x}_2\|^2_{{\bs \Sigma}^{-1}_i} - \Delta^2_{{\bs \Sigma}^{-1}_i}\right| &\xrightarrow{\mathcal{P}}& 0,
\end{eqnarray*}
so that by Slutsky's lemma \citep[][p 11]{vdv98} it follows that, as $n_i, p \rightarrow \infty$,
\[
\underline{\epsilon}_n \xrightarrow{\mathcal{P}} \frac{1}{2}\left[\Phi\left(-\frac{\Delta^2_{0, {\bf I}}}{\Delta_{{\bs \Sigma}^{-1}_1}}\right) + \Phi\left(\frac{\Delta^2_{0, {\bf I}}}{\Delta_{{\bs \Sigma}^{-1}_2}}\right)\right],
\]
where the convergence remains (asymptotically) true even for the sample based classifier $\underline{A}_0({\bf x})$ since $\epsilon_o^{\text{ideal}}$ is the limiting value of $\underline{\epsilon}_n$.\\
%%%%%%%%%%%%%%%%%%%%%%%%%%%%%%%%%%%%%%%%%%%%%%%%%%%%%%%%%%%%%%%
\indent Finally, we consider a similar evaluation of the classifier under non-normality which we discuss for the general class of elliptically contoured distributions (including multivariate normal). Let $\mathcal{F}_{i}\equiv \mathcal{E}_{p} (h, {\bs \mu}_{i}, {\bs \Sigma}_{i})$, $i=1,2$,
where $\mathcal{E}_{p} (h, {\bs \mu}, {\bs \Sigma})$ denotes a $p$-dimensional elliptical distribution with density function
\begin{eqnarray}
|\bs \Sigma|^{-1/2} h\left ( \| {\bf x} - {\bs \mu} \|_{{\bs \Sigma}}^{2} \right),\label{eqn: ElDens}
\end{eqnarray}
where $h$ is a monotone (decreasing) function on $[0,\infty)$ and parameters ${\bs \mu}$ and ${\bs \Sigma}$ are specified as in model (7). Assume now that ${\bf x}_{ik} \sim \mathcal{E}_{p} (h,{\bs \mu}_{i}, {\bs \Sigma_{i}})$ and denote the first two conditional moments of the classifier $A_{0}(\bs x)$ by
\[
E_{i}= \text{E} \left[ A_{0}({\bf x})| \overline{\bf x}_i, U_{n_i} \right ], \quad V_{i}= \text{Var} \left [A_{0}({\bf x})| \overline{\bf x}_i, U_{n_i} \right ],
\]
respectively, $i = 1, 2$. \citet{Wakaki94} has discussed Fisher's linear discriminant function for elliptical distributions and also its robustness particularly in the context of $M$-estimation. Our assertion in the following closely follows his structure. \\
%%%%%%%%%%%%%%%%%%%%%%%%%%%%%%%%%%%%%%%%%%%%%%%%%%%%%%%%%%%%%%%
\indent By conditioning the classifier on the data, i.e. by considering $\underline{A}_0({\bf x}) = A_0({\bf x})|(\overline{\bf x}_i, U_{n_i})$ and using its standardized version, we get \citep[see][Theorem 1.1, p 260]{Wakaki94}
\begin{eqnarray}
P \left ( \frac{\underline{A}_0({\bf x}) -E_{i}}{\sqrt{V_{i}}} \leq z \right ) =: Q(z) \label{eqn: Qdistr}
\end{eqnarray}
for any $z \in R$, where $\Gamma$ denotes the gamma function, $h$ is defined in (\ref{eqn: ElDens}), and $Q$ is the distribution function whose density function is given by
\[
q(z)= \frac{\pi^{(p - 1)/2}}{\Gamma\left(\frac{p-1}{2}\right)} \int_{0}^{\infty} s^{(p-3)/2} h(z^2 + s)ds.
\]
The normal distribution is a special case of $\mathcal{E}_{p}$, for if $h(s)=(2\pi)^{-p/2} \exp(-s/2)$, then $q(z)$ reduces to the standard normal. Now, using (\ref{eqn: Qdistr}), we can express the conditional (or actual) misclassification probability as
 \begin{eqnarray*}
\underline{\epsilon}_n = \frac{1}{2}\left[Q\left(-\frac{E_{1}}{\sqrt{V_{1}}}\right) + Q \left(\frac{E_{2}}{\sqrt{V_{2}}}\right)\right].
\end{eqnarray*}

%%%%%%%%%%%%%%%%%%%%%%%%%%%%%%%%%%%%%%%%%%%%%%%%%%%%%%%%%%%%%%%
\section{Multi-sample $U$-classifier}\label{sec:MultSampDF}
%%%%%%%%%%%%%%%%%%%%%%%%%%%%%%%%%%%%%%%%%%%%%%%%%%%%%%%%%%%%%%%

It is obvious from the construction of the two-sample classifier that it can be easily extended to the multi-sample case. Let
\[
{\bf x}_{ik} \sim \mathcal{F}_i, ~k = 1, \ldots, n_i
\]
be iid vectors with $\text{E}({\bf x}_{ik}) = {\bs \mu}_i$, $\text{Cov}({\bf x}_{ik}) = {\bs \Sigma}_i$, $i = 1, \ldots, g \geq 2$. The multi-sample version of classifier in Eqn. (\ref{eqn:ModifiedTwoSDF}) can be expressed as
\begin{equation}\label{eqn:MultSDFDirectExt}
A_0({\bf x}) = \frac{1}{p}{\bf x}'_0 (\overline{\bf x}_i - \overline{\bf x}_k) - \frac{1}{2}(U_{n_i} - U_{n_k}),
\end{equation}
$i, k = 1, \ldots, g$, $k \neq i$. Alternatively, to write it in a more explicit form, let
\[
A_{0i}({\bf x}) = \frac{1}{p}{\bf x}'_0\overline{\bf x}_i - \frac{1}{2}U_{n_i}
\]
be the discriminant function for population $i$, so that the classifier is
\begin{equation}\label{eqn:MultSDFSimplerForm}
A_{0g}({\bf x}) = A_{0i}({\bf x}) - A_{0k}({\bf x})
\end{equation}
for any distinct pair $(i, k)$. The classification rule modifies to
%%%%%%%%%%%%%%%%%%%%%%%%%%%%%%%%%%%%%%%%%%%%%%%%%%%%%%%%%%%%%%%
\[
\text{Assign ${\bf x}_0$ to $\pi_i$ if $A_{0g}({\bf x}) > 0$, i.e., if $A_{0i}({\bf x}) > A_{0k}({\bf x})$; otherwise to $\pi_k$.}
\]
%%%%%%%%%%%%%%%%%%%%%%%%%%%%%%%%%%%%%%%%%%%%%%%%%%%%%%%%%%%%%%%
To study the properties of the multi-sample classifier and its asymptotic behavior, we first extend the two-sample assumptions for the general case.
%------------------------------------------------------------%
\begin{assn}\label{assn:4thMomntMS}
$\text{E}(x^4_{iks}) = \gamma < \infty$, $\gamma \in \mathbb{R}^+$, $i = 1, \ldots, g$.
\end{assn}
\begin{assn} \label{assn:traceSigmaMS}
$\lim_{p \rightarrow \infty}\text{tr}({\bs \Sigma}_i)/p = O(1)$, $i = 1, \ldots, g$.
\end{assn}
\begin{assn}\label{assn:ExtraH1DistnMS}
$\lim_{p \rightarrow \infty}{\bs \mu}'_i{\bs \Sigma}_k{\bs \mu}_j/p = O(1)$, $i, j, k = 1, \ldots, g$.
\end{assn}
\begin{assn}\label{assn:traceSigmaHadProdMS}
$\lim_{p \rightarrow \infty}\frac{\text{tr}({\bs \Sigma}^a_i\odot{\bs \Sigma}^b_j)}{\text{tr}({\bs \Sigma}^a_i\otimes{\bs \Sigma}^b_j)} = 0$, $a, b = 1, 2, 3$, $a + b \leq 4$, $i, j = 1, \ldots, g$, where $\odot$ and $\otimes$ are Hadamard and Kronecker products. %$i = j$ or $i \neq j$,
\end{assn}
%--------------------------------------------------------------%
We begin by the following lemma which generalizes Lemma \ref{lem:MomsTwoSDF} on the moments of the two-sample classifier.
%%%%%%%%%%%%%%%%%%%%%%%%%%%%%%%%%%%%%%%%%%%%%%%%%%%%%%%%%%%%%%%
%%%%%%%%%%%%%%%%%%%%%%%%%%%%%%%%%%%%%%%%%%%%%%%%%%%%%%%%%%%%%%%
\begin{lem}\label{lem:MomsMultSDF} Given $A_0({\bf x})$ in Eqn. (\ref{eqn:MultSDFDirectExt}) or (\ref{eqn:MultSDFSimplerForm}). Let ${\bf x}_0 \in \pi_i$. Then
\begin{eqnarray}
\text{E}[A_0({\bf x})|\pi_i] &=& \|{\bs \mu}_i - {\bs \mu}_k\|^2/2p = \Delta^2_{\bf I}/2p\label{eqn:MeanDFMultS}\\
\text{Var}[A_0({\bf x})|\pi_i] &=& \delta^2_i/p^2 + \|{\bs \mu}_i - {\bs \mu}_k\|^2_{{\bs \Sigma}^{-1}_i}/p^2 = \delta^2_i/p^2 + \Delta^2_{{\bs \Sigma}^{-1}_i}/p^2,\label{eqn:VarDFMultS}
\end{eqnarray}
where $\Delta^2_{{\bf M}^{-1}} = ({\bs \mu}_i - {\bs \mu}_k)'{\bf M}({\bs \mu}_i - {\bs \mu}_k)$, ${\bf M} = {\bs \Sigma}$ or ${\bf M} = {\bf I}$, and
\[
\delta^2_i = \frac{\text{tr}({\bs \Sigma}^2_i)}{n_i} + \frac{\text{tr}({\bs \Sigma}_i{\bs \Sigma}_k)}{n_k} + \sum_{i = 1}^2\frac{\text{tr}({\bs \Sigma}^2_i)}{2n_i(n_i - 1)}, ~~i, k = 1, \ldots, g, ~i \neq k.
\]
\end{lem}
%%%%%%%%%%%%%%%%%%%%%%%%%%%%%%%%%%%%%%%%%%%%%%%%%%%%%%%%%%%%%%%
%%%%%%%%%%%%%%%%%%%%%%%%%%%%%%%%%%%%%%%%%%%%%%%%%%%%%%%%%%%%%%%
The moment estimators follow obviously from those of the two-sample case given in Eqns. (\ref{eqn:LocInvE0})-(\ref{eqn:LocInvE2}). Likewise, the consistency of these estimators follows from Lemma \ref{lem:MomsTwoSDF}. This helps us extend Theorem \ref{thm:ConsistencyTwoSDF} on the consistency and asymptotic normality of $A_0({\bf x})$ for the general case as following.
%%%%%%%%%%%%%%%%%%%%%%%%%%%%%%%%%%%%%%%%%%%%%%%%%%%%%%%%%%%%%%%
%%%%%%%%%%%%%%%%%%%%%%%%%%%%%%%%%%%%%%%%%%%%%%%%%%%%%%%%%%%%%%%
\begin{thm}\label{thm:ConsisAsympNDMultSDF} Given $A_0({\bf x})$ in Eqn. (\ref{eqn:MultSDFDirectExt}) or (\ref{eqn:MultSDFSimplerForm}) with its moments in Lemma \ref{lem:MomsTwoSDF}. Let ${\bf x}_0 \in \pi_i$. Under Assumptions \ref{assn:4thMomntMS}-\ref{assn:ExtraH1DistnMS}, as $n_i, p \rightarrow \infty$, $i = 1, \ldots, g$,
\begin{eqnarray}
\frac{A_0({\bf x})}{\Delta^2_{\bf I}/p} &\xrightarrow{\mathcal{P}}& \frac{(-1)^i}{2} + o_P(1),\\
\frac{A_0({\bf x}) - \text{E}[A_0({\bf x})]}{\sqrt{\text{Var}[A_0({\bf x})]}} &\xrightarrow{\mathcal{D}}& N(0, 1),
\end{eqnarray}
Further, the limits hold when the moments are replaced with their empirical estimators given in Eqns. (\ref{eqn:LmtgMeanDFEstd})-(\ref{eqn:LmtgVarDFEstd}).
\end{thm}
%%%%%%%%%%%%%%%%%%%%%%%%%%%%%%%%%%%%%%%%%%%%%%%%%%%%%%%%%%%%%%%
%%%%%%%%%%%%%%%%%%%%%%%%%%%%%%%%%%%%%%%%%%%%%%%%%%%%%%%%%%%%%%%
As the multi-sample case is a straightforward extension of its two-sample counterpart in Sec. \ref{sec:TwoSampDF}, we skip many detailed proofs to avoid unnecessary repetitions.

%%%%%%%%%%%%%%%%%%%%%%%%%%%%%%%%%%%%%%%%%%%%%%%%%%%%%%%%%%%%%%
\section{Simulations}\label{sec:Simns}
%%%%%%%%%%%%%%%%%%%%%%%%%%%%%%%%%%%%%%%%%%%%%%%%%%%%%%%%%%%%%%

We use simulation results to evaluate the performance of $A_0({\bf x})$ under practical scenarios, mainly focusing on consistency, asymptotic normality and control of misclassification under high-dimensional framework. We consider $g = 2$ case and generate data from multivariate normal and $t$ distributions, i.e. $\mathcal{F}_i$ is either $\mathcal{N}_p({\bs \mu}_i, {\bs \Sigma}_i)$, $i = 1, 2$, or $t_\nu({\bs \mu}_i, {\bs \Sigma}_i)$, $\nu = 10$, $i = 1, 2$. For each distribution, we set ${\bs \mu}_1 = {\bf 0}$ with $\lfloor p/3\rfloor$ elements of ${\bs \mu}_2$ also 0 and the rest as 1 where $\lfloor\cdot\rfloor$ denotes the smallest integer. For ${\bs \Sigma}_i$, we consider two cases: (1) Both populations have AR(1) structure, $\text{Cov}(X_k, X_l) = \kappa\rho^{|k - l|}$, $\forall~ k, l$, with $\sigma^2 = 1$ for $i = 1$ and 2, where $\rho = 0.3$ for $i = 1$ and $\rho = 0.7$ for $i = 2$, to represent both low and high correlation structures; (2) The same AR(1) structure for $i = 1$ with $\sigma^2 = 1$, $\rho = 0.5$, whereas an unstructured (UN) ${\bs \Sigma}_i$ for $i = 2$, defined as ${\bs \Sigma} = (\sigma_{ij})_{i, j = 1}^p$ with $\sigma_{ii} = 1(1)p$ and $\rho_{ij} = (i - 1)/p$, $i \neq j$. \\
%%%%%%%%%%%%%%%%%%%%%%%%%%%%%%%%%%%%%%%%%%%%%%%%%%%%%%
\indent For finite-sample performance of the classifier under arbitrarily growing dimension, emphasizing the $p \gg n_i$, we generate samples of sizes $n_1 = 5$, $n_2 = 7$, and combine with $p = \{10, 20, 50, 100, 300, 500, 700, 1000, 3000, 5000, 10000\}$. Finally, all results are averages of 1000 simulation runs for each combination of parameters mentioned above. To additionally observe the effect of large $n_i$, the misclassification rates are also presented for $n_1 = 10, n_2 = 12$. Similarly, to assess the classifier for very different sample sizes, we also used $n_1 = 5, n_2 = 25$ and $n_1 = 10, n_2 = 50$ and observed very similar results, hence not reported here. \\
%%%%%%%%%%%%%%%%%%%%%%%%%%%%%%%%%%%%%%%%%%%%%%%%%%%%%%%%%%%%%%%%%%%%%%%
\indent Fig. 2 shows the results of asymptotic normality of $A_0({\bf x})$, where the first two rows are for normal distribution, respectively for AR-AR and AR-UN structures. Each row gives a histogram of $A_0({\bf x})$, with empirical density added to it, for $p = 100, 500$ and 1000 (left to right). The last two rows are for multivariate $t$ distribution with $\nu = 10$. As stated above, the results are carried out for several other dimensions as well, up to $p = 10000$, and also for other sample sizes, but due to similarity of the graphs, only a selection is reported here.\\
%%%%%%%%%%%%%%%%%%%%%%%%%%%%%%%%%%%%%%%%%%%%%%%%%%%%%%%%%%%%%%%%%%%%%%%
\indent We observe close normal approximation for $n_i$ as small as 5 or 7, and the results for $t$ distribution depict small sample robustness of the classifier to non-normality. To make the results of the two distributions comparable, the density axes are scaled to the same height, so that the heavy-tailed behavior of $t$ distribution can be witnessed from a slightly extended range on the x-axis. In general, it is observed that the increasing dimension does not damage the asymptotic normality of the classifier even if the data are non-normal.\\
%%%%%%%%%%%%%%%%%%%%%%%%%%%%%%%%%%%%%%%%%%%%%%%%%%%%%%%%%%%%%%
\indent A similar performance is observed for the control of misclassification rate, shown in Figs. 3 for $n_1 = 5$, $n_2 = 7$, and 4 for $n_1 = 10, n_2 = 12$. The thick line represents the actual error rate under asymptotic normality of the classifier, i.e. $\Phi(-\text{E}(A_0)/\sqrt{\text{Var}(A_0)})$, where $\Phi(\cdot)$ is the (univariate) normal distribution function and $\text{E}(\cdot)$, $\text{Var}(\cdot)$ are the moments of the classifier. This actual error rate is used as a reference to assess the estimated error rate shown in dashed (dotted) line for normal ($t_{10}$) distribution. Further, the upper (lower) panel in each figure is for AR-AR (AR-UN) pair of covariances.\\
%%%%%%%%%%%%%%%%%%%%%%%%%%%%%%%%%%%%%%%%%%%%%%%%%%%%%%%%%%%%%%
\indent The estimated error closely follows the actual error for $n_1 = 5, n_2 = 7$, and the error rate also converges to zero, showing consistency of the classifier. For $t$
%%%%%%%%%%%%%%%%%%%%%%%%%%%%%%%%%%%%%%%%%%%%%%%%%%%%%%%%%%%%%%
\begin{landscape}
\begin{figure}[h!]\centering\label{fig:FigAsympND}
\includegraphics[width = 0.55\textwidth]{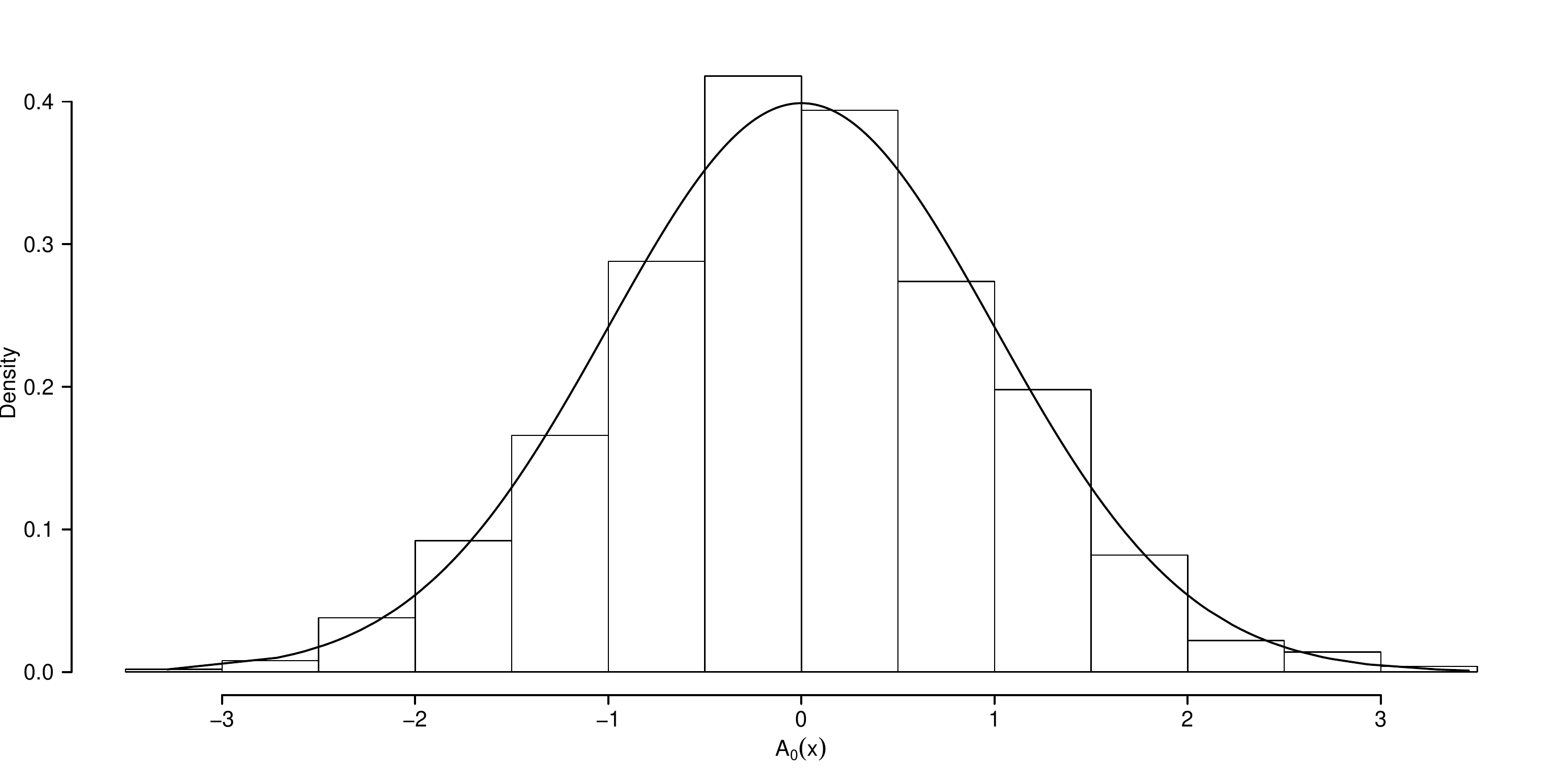}~~
\includegraphics[width = 0.55\textwidth]{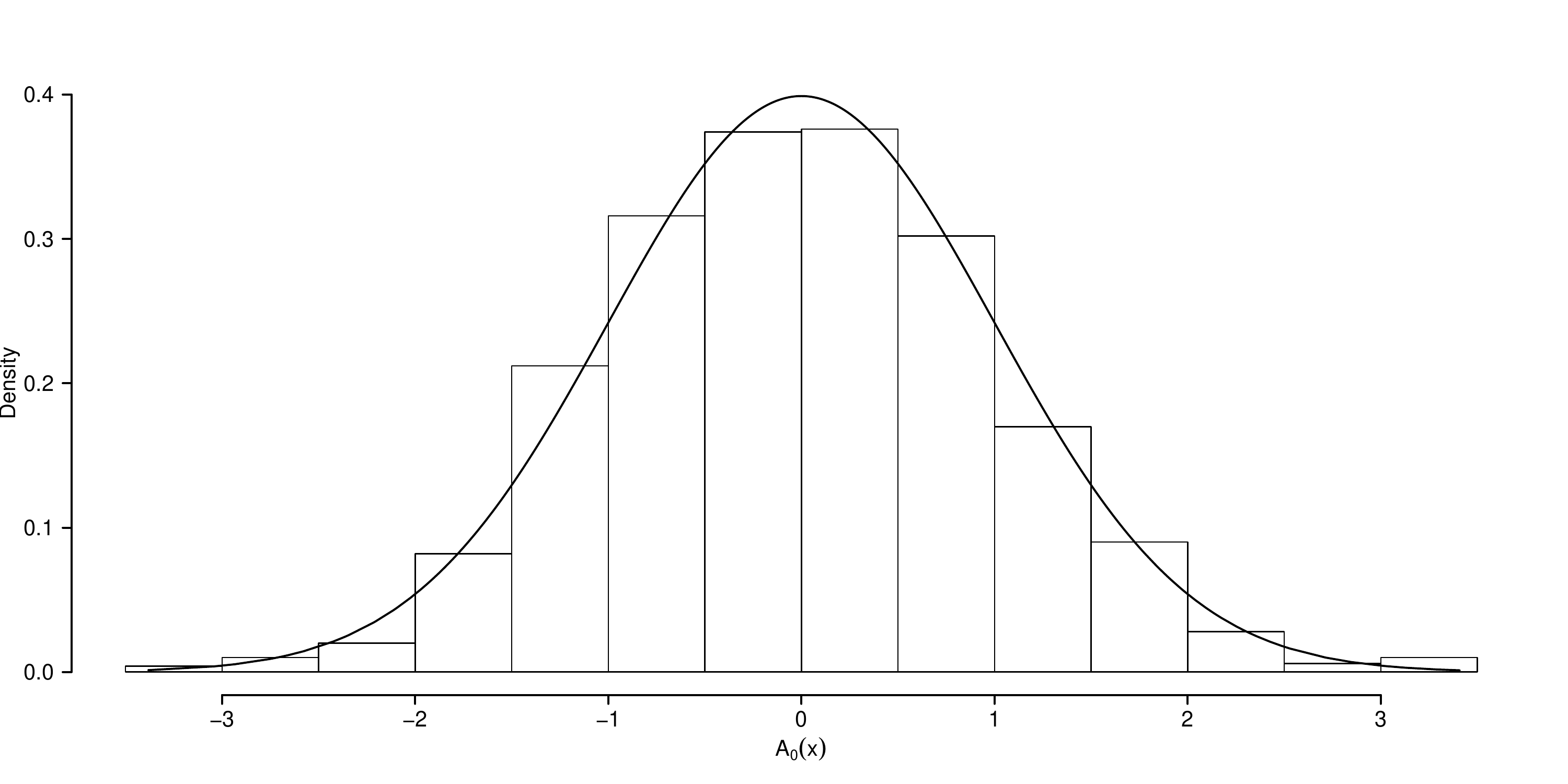}~~
\includegraphics[width = 0.55\textwidth]{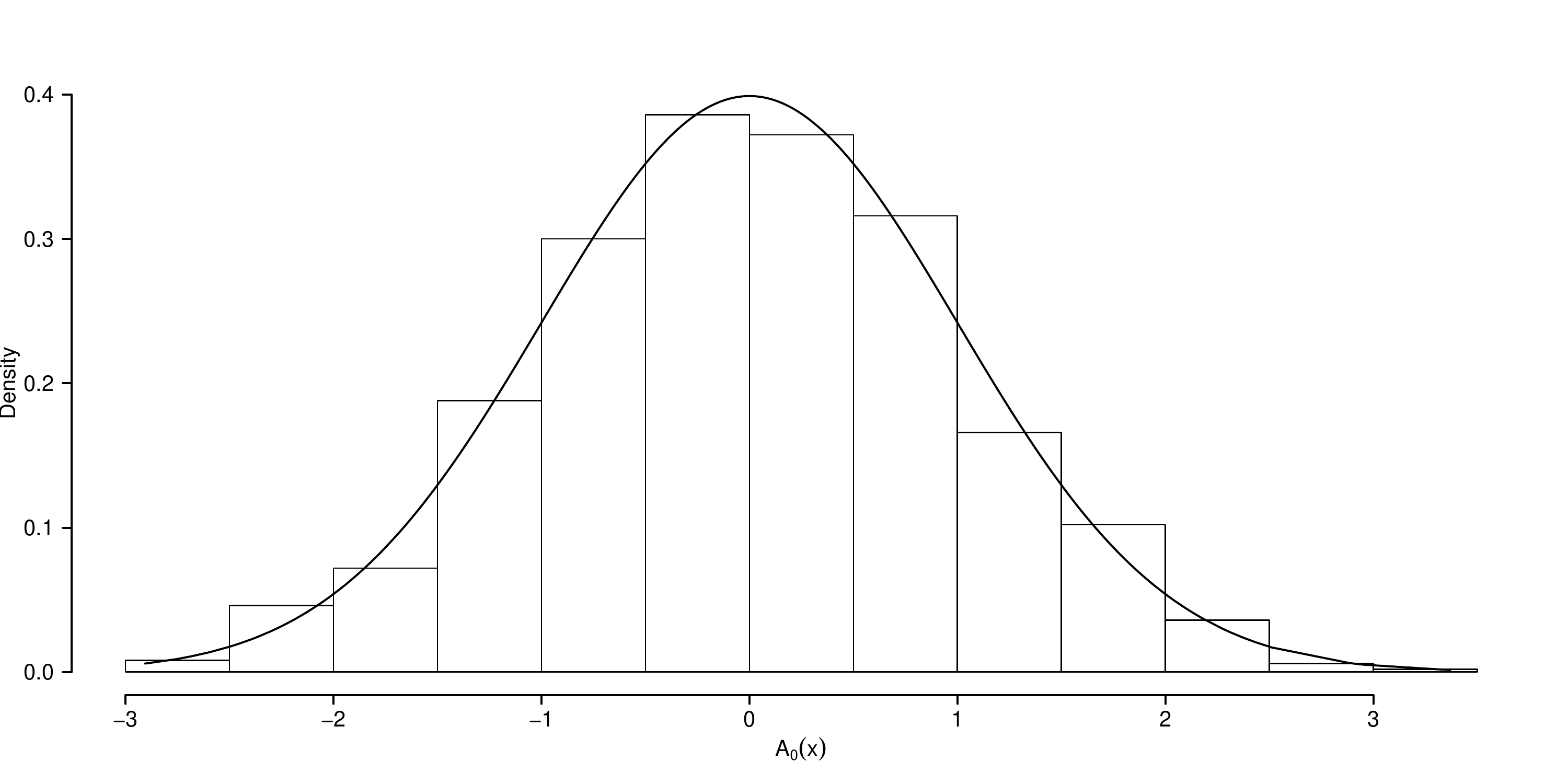}\\
\includegraphics[width = 0.55\textwidth]{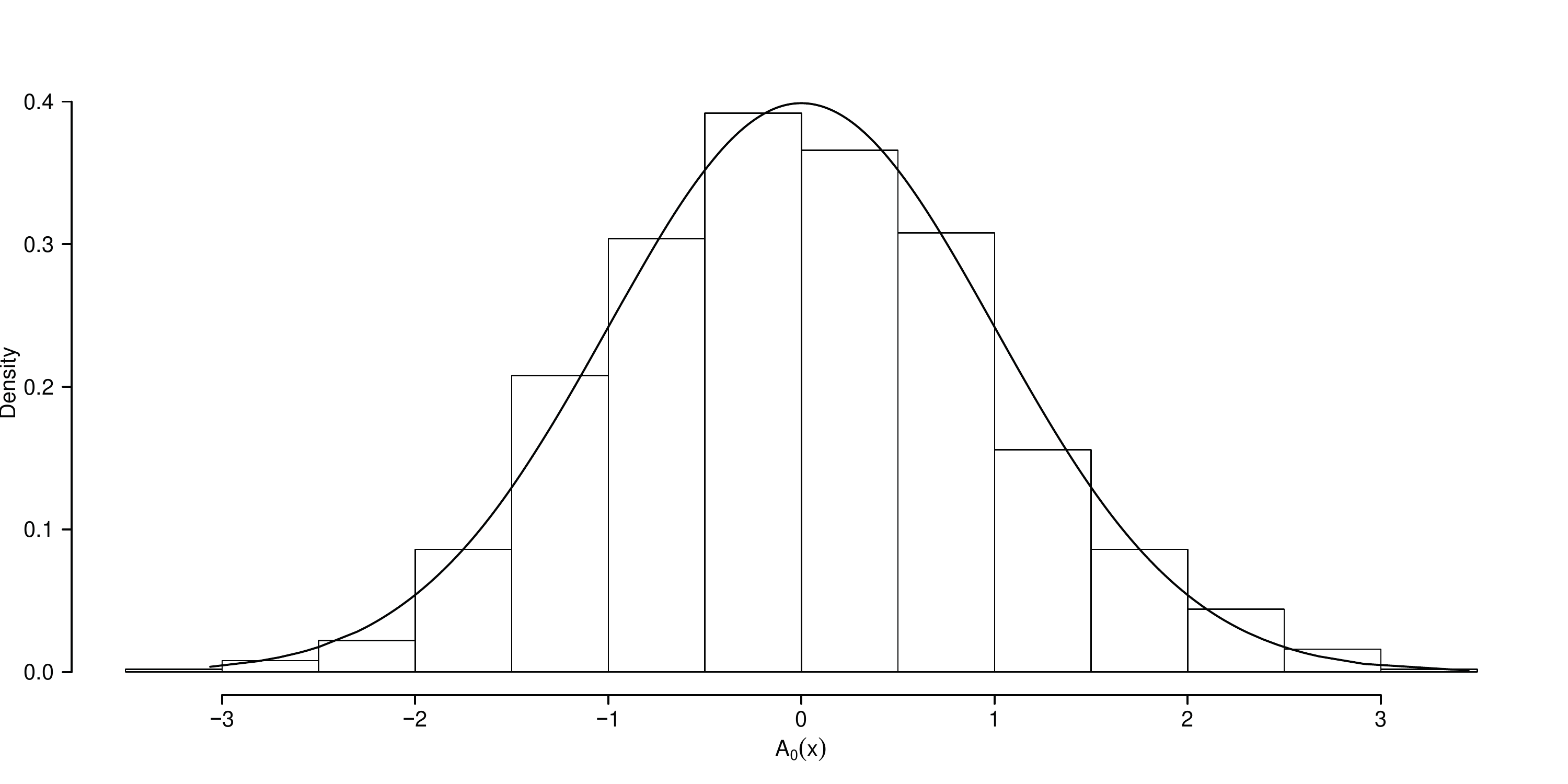}~~
\includegraphics[width = 0.55\textwidth]{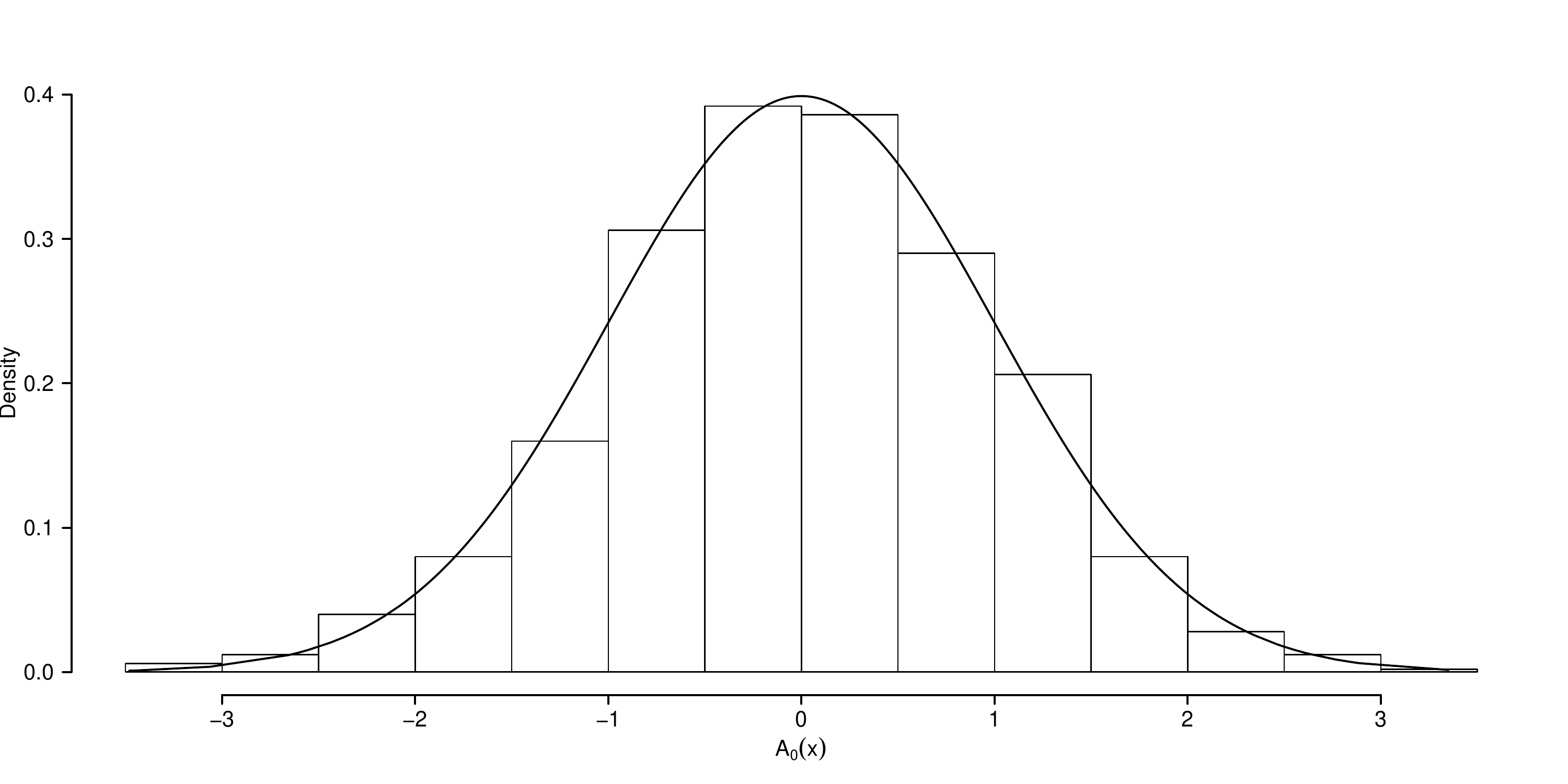}~~
\includegraphics[width = 0.55\textwidth]{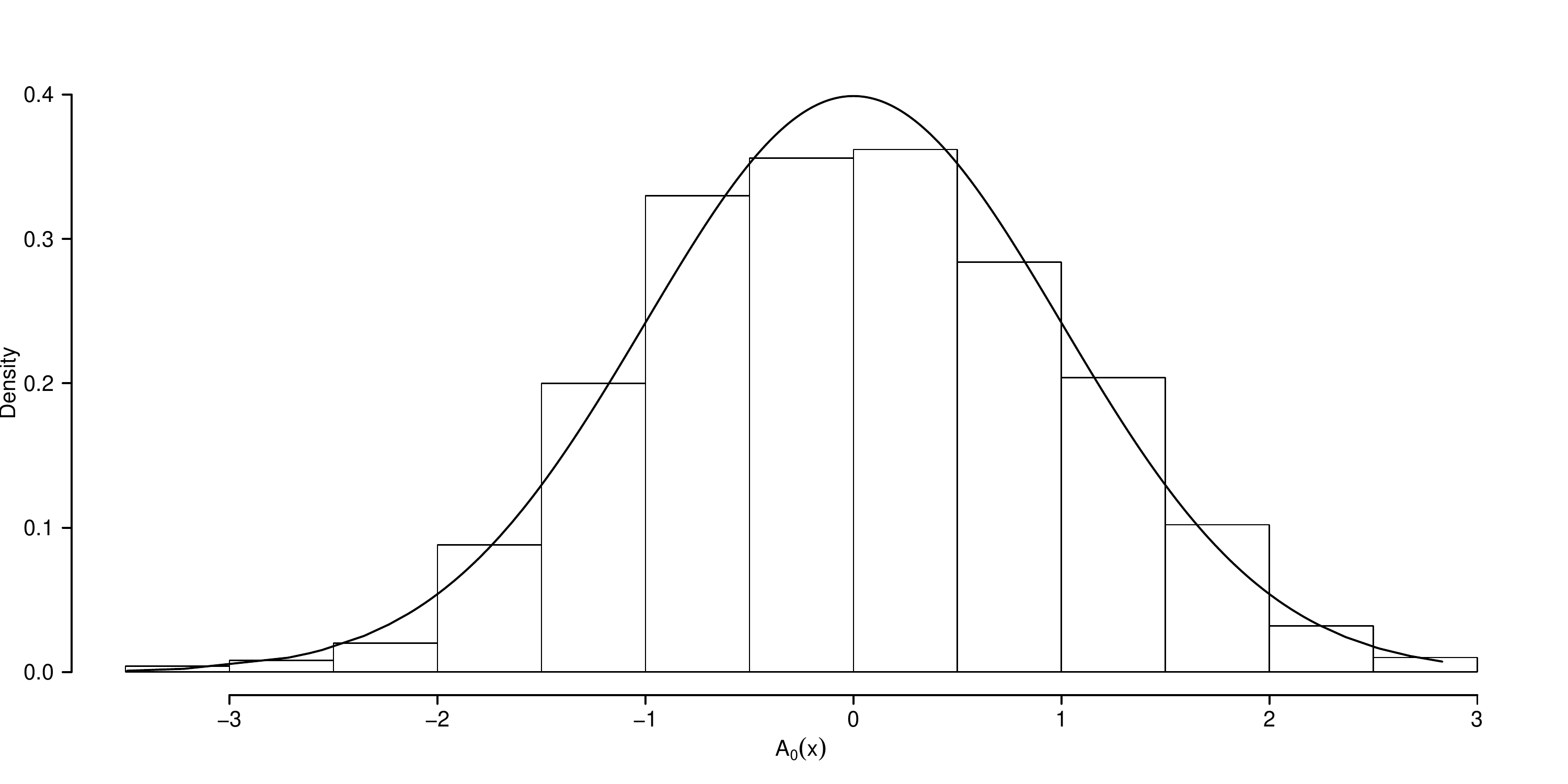}
\includegraphics[width = 0.55\textwidth]{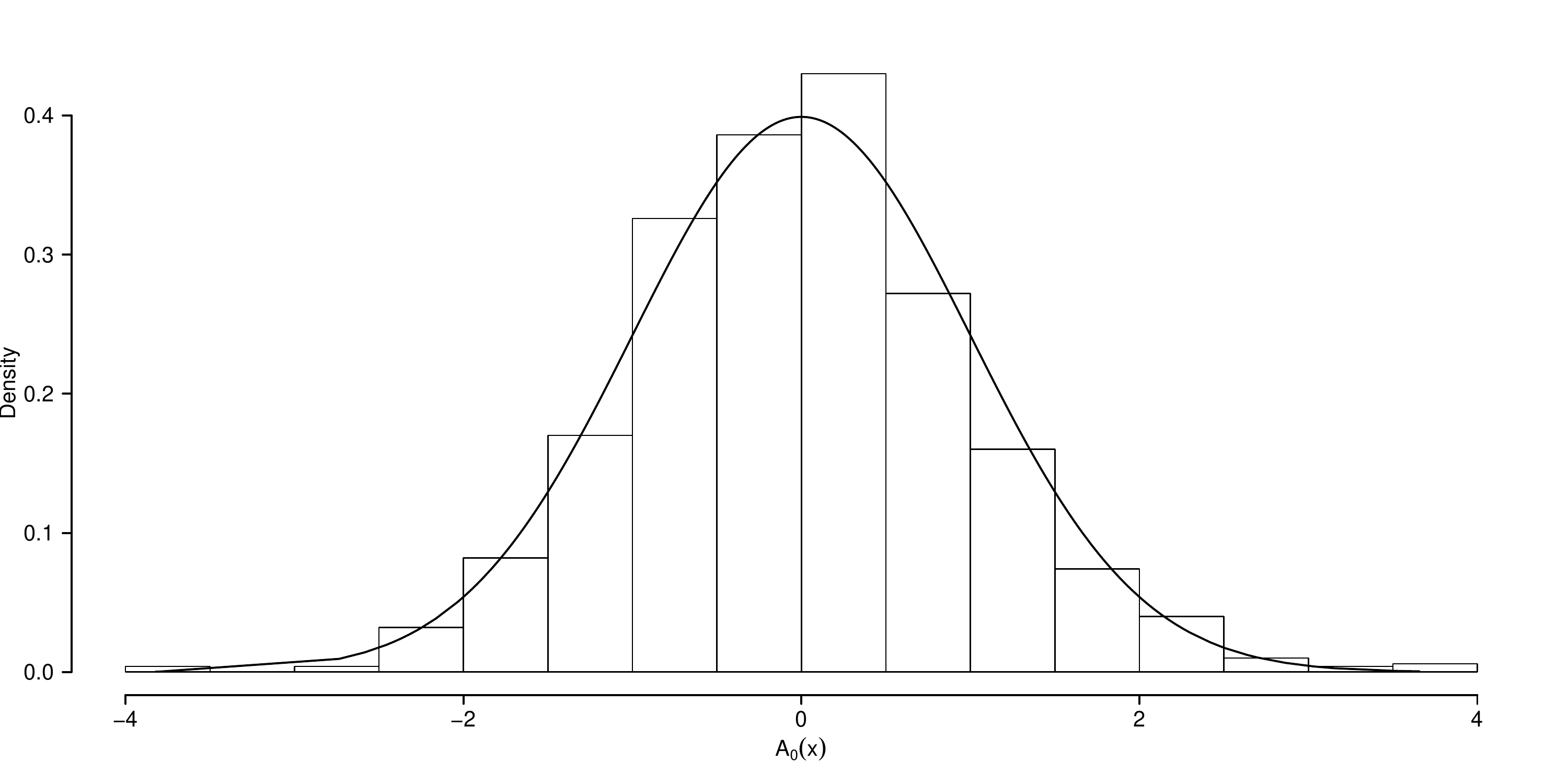}~~
\includegraphics[width = 0.55\textwidth]{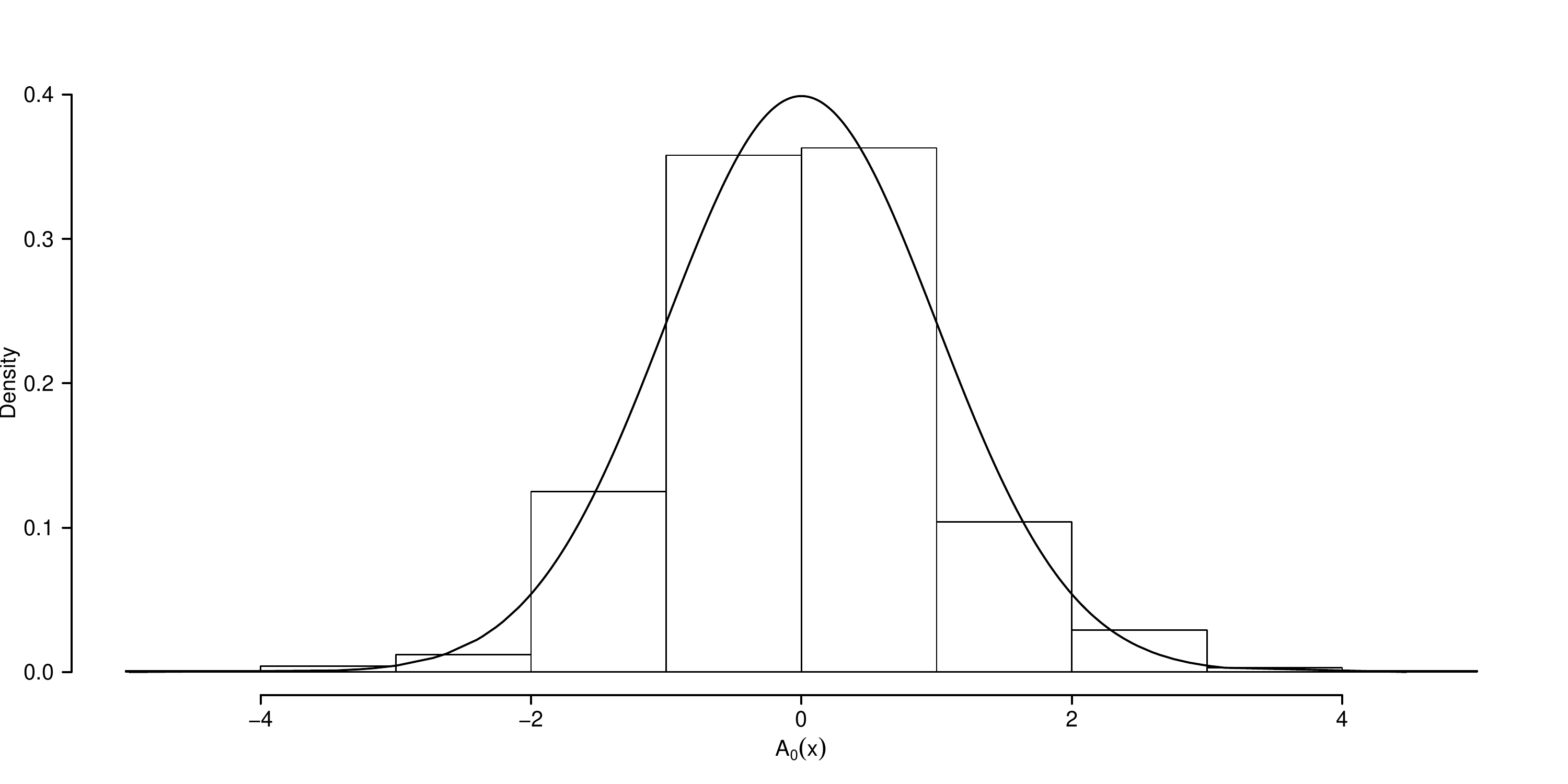}~~
\includegraphics[width = 0.55\textwidth]{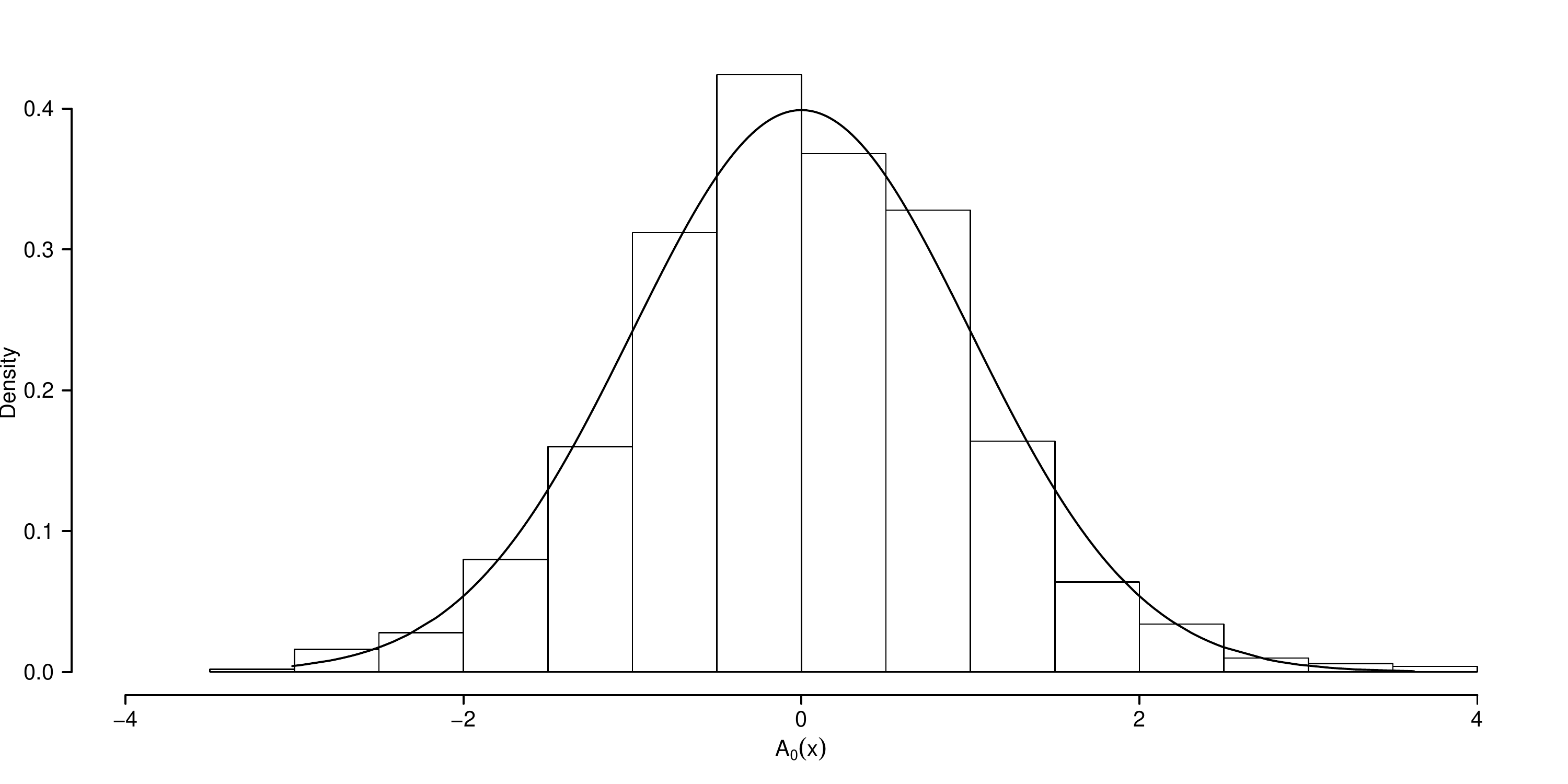}\\
\includegraphics[width = 0.55\textwidth]{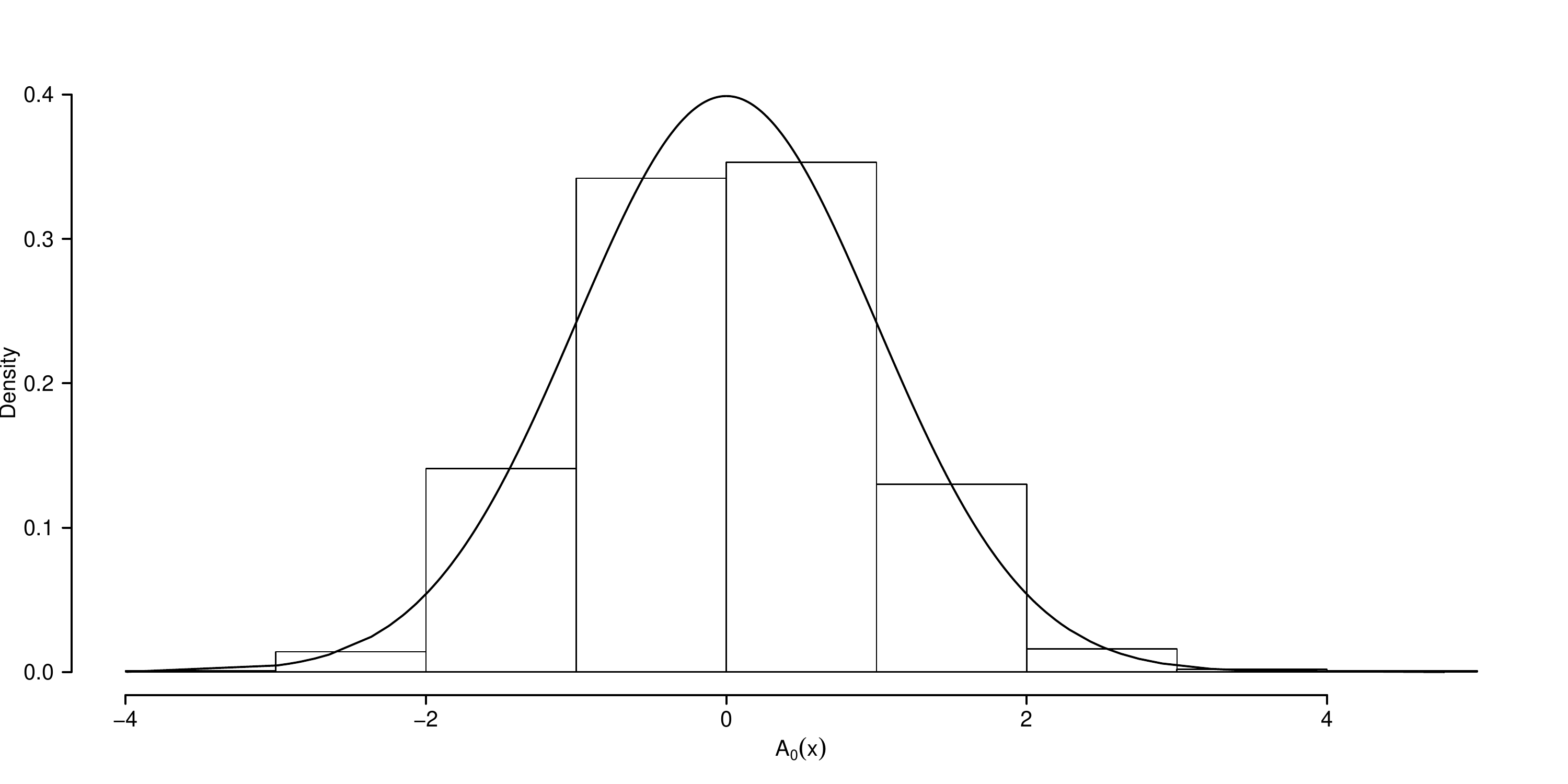}~~
\includegraphics[width = 0.55\textwidth]{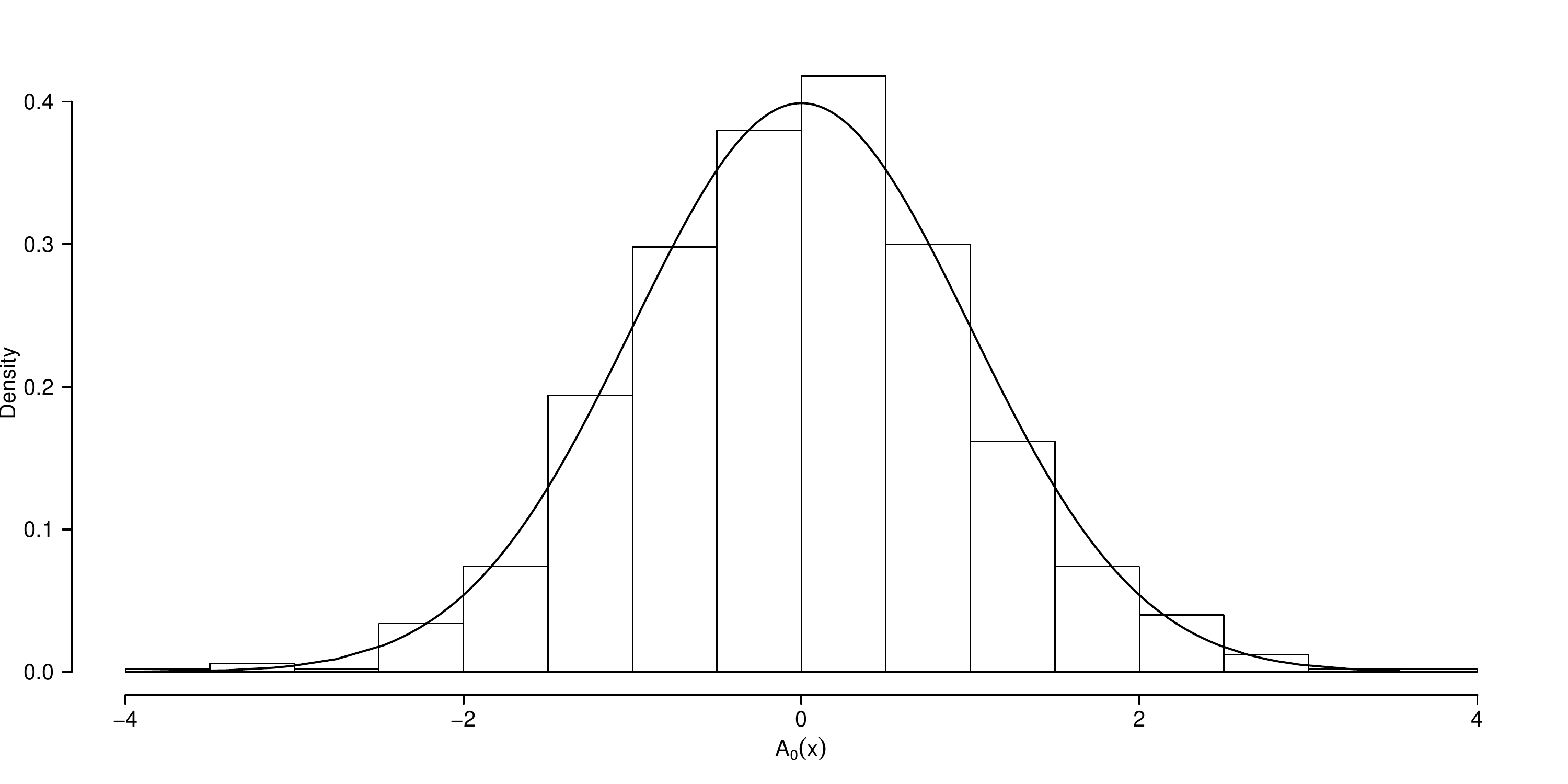}~~
\includegraphics[width = 0.55\textwidth]{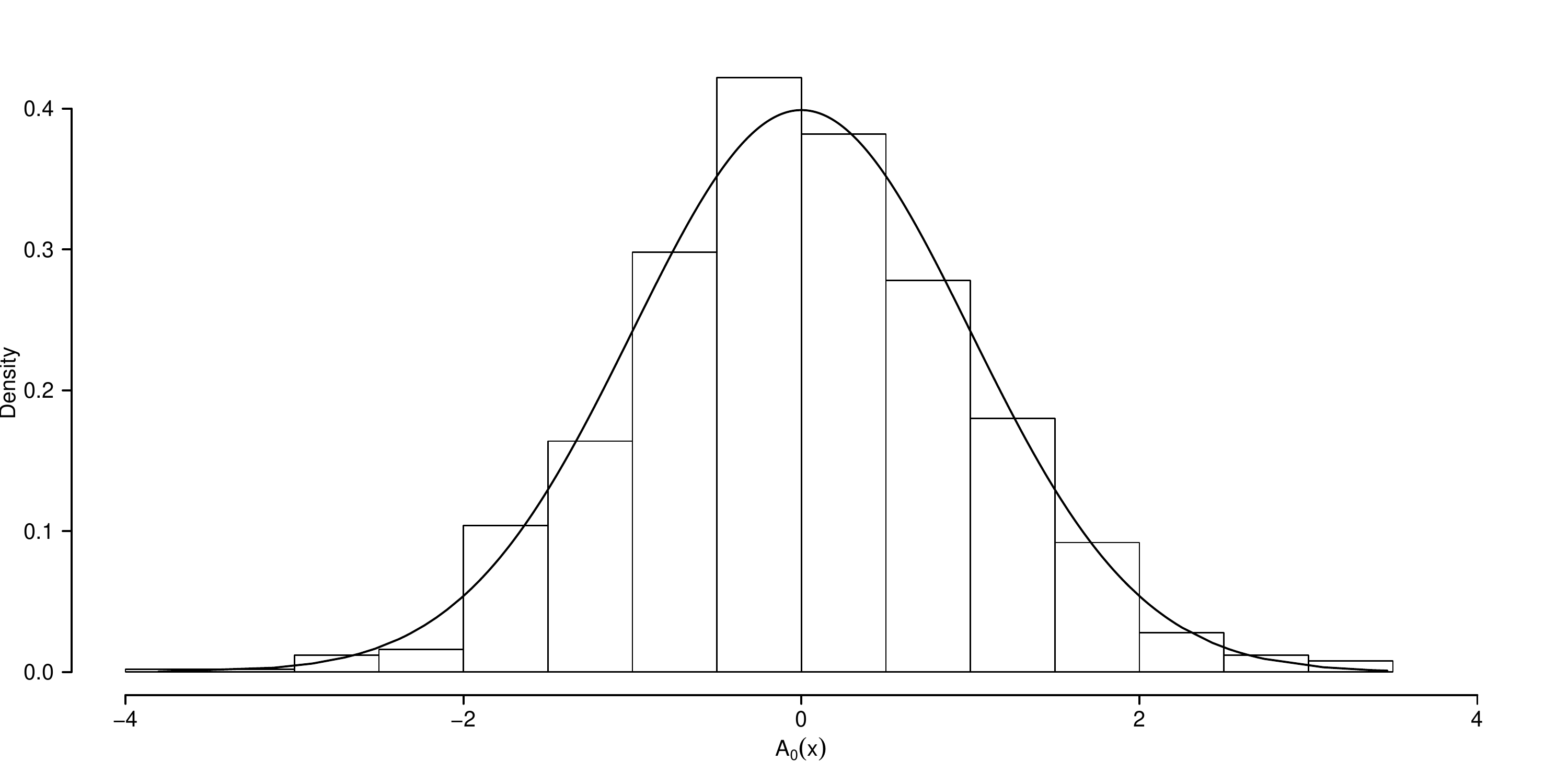}
\caption{\footnotesize \sc Asymptotic normality of $A_0$ for two-class case with data from multivariate normal (Rows 1-2) \& $t$ (Rows 3-4) distributions, $n_1 = 5$, $n_2 = 7$, $p = 100, 500, 1000$ (L to R each row) and covariance structures AR-AR (Rows 1, 3) \& AR-UN (Rows 2, 4).}
\end{figure}
\end{landscape}
%%%%%%%%%%%%%%%%%%%%%%%%%%%%%%%%%%%%%%%%%%%%%%%%%%%%%%%%%%%%%%%%%%%%%
\begin{figure}[t!]\centering\label{fig:FigAPERsmallN}
\includegraphics[width = 1.025\textwidth]{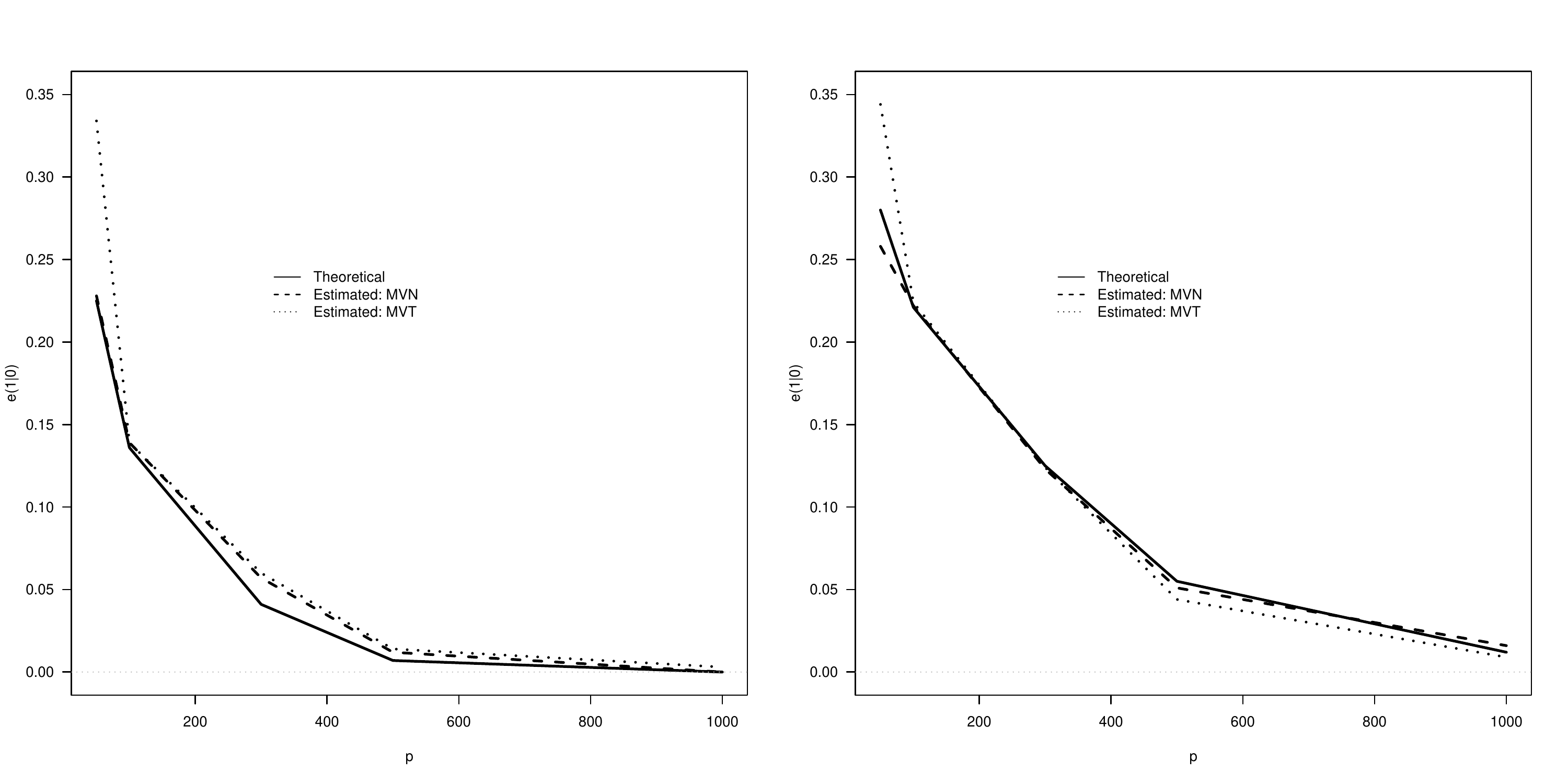}
\includegraphics[width = 1.025\textwidth]{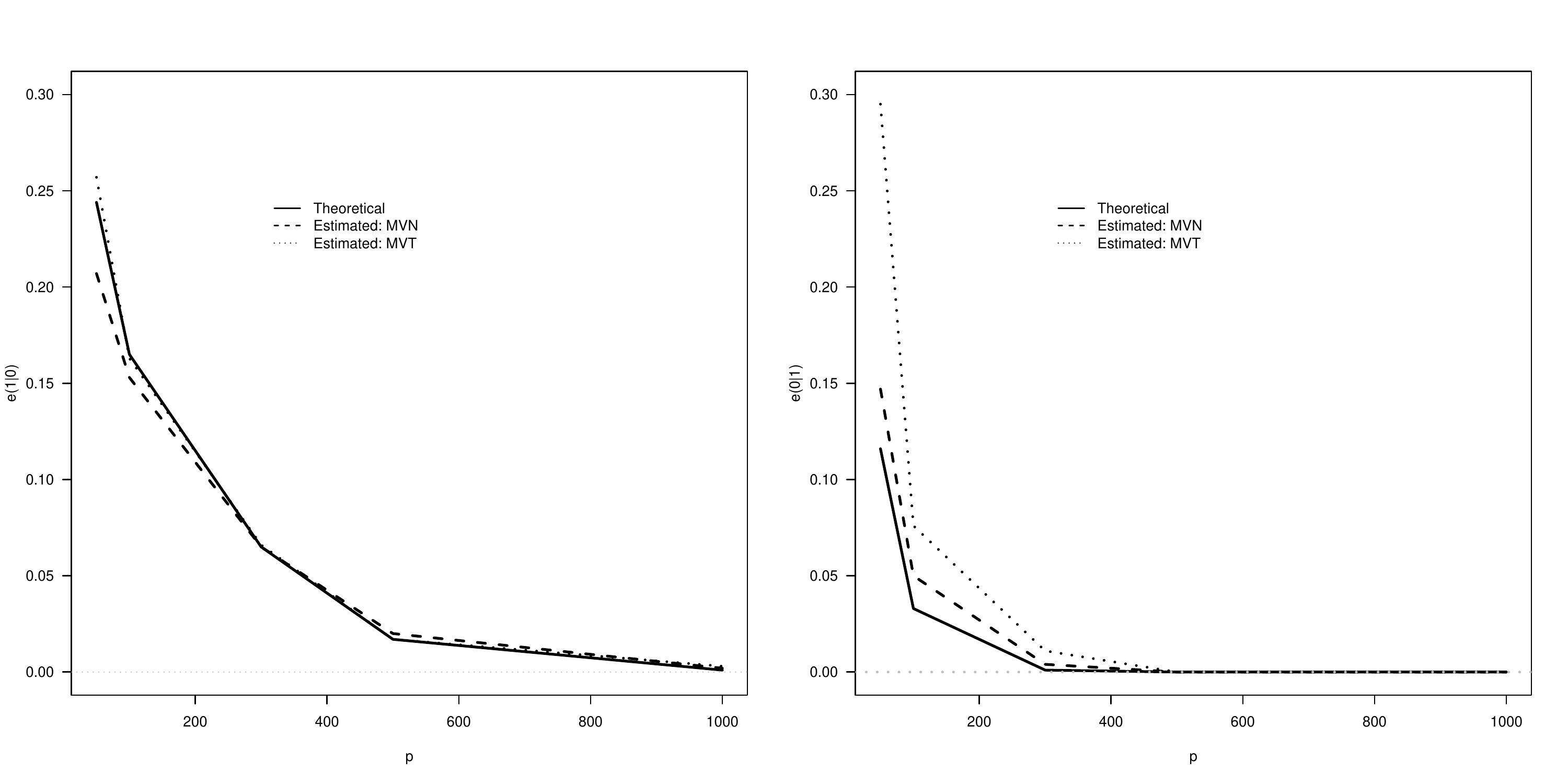}
\caption{\footnotesize \sc Theoretical (thick line) and Estimated error rates of $A_0({\bf x})$ for two-class case with data from multivariate normal (dashed line) and $t$ (dotted line) distributions, $n_1 = 5$, $n_2 = 7$, $p = 10, 20, 50, 100, 200, 300, 500$ and covariance structures AR-AR (upper panel) \& AR-UN (lower penal).}
\end{figure}
\begin{figure}[t!]\centering\label{fig:FigAPERlargeN}
\includegraphics[width = 1.025\textwidth]{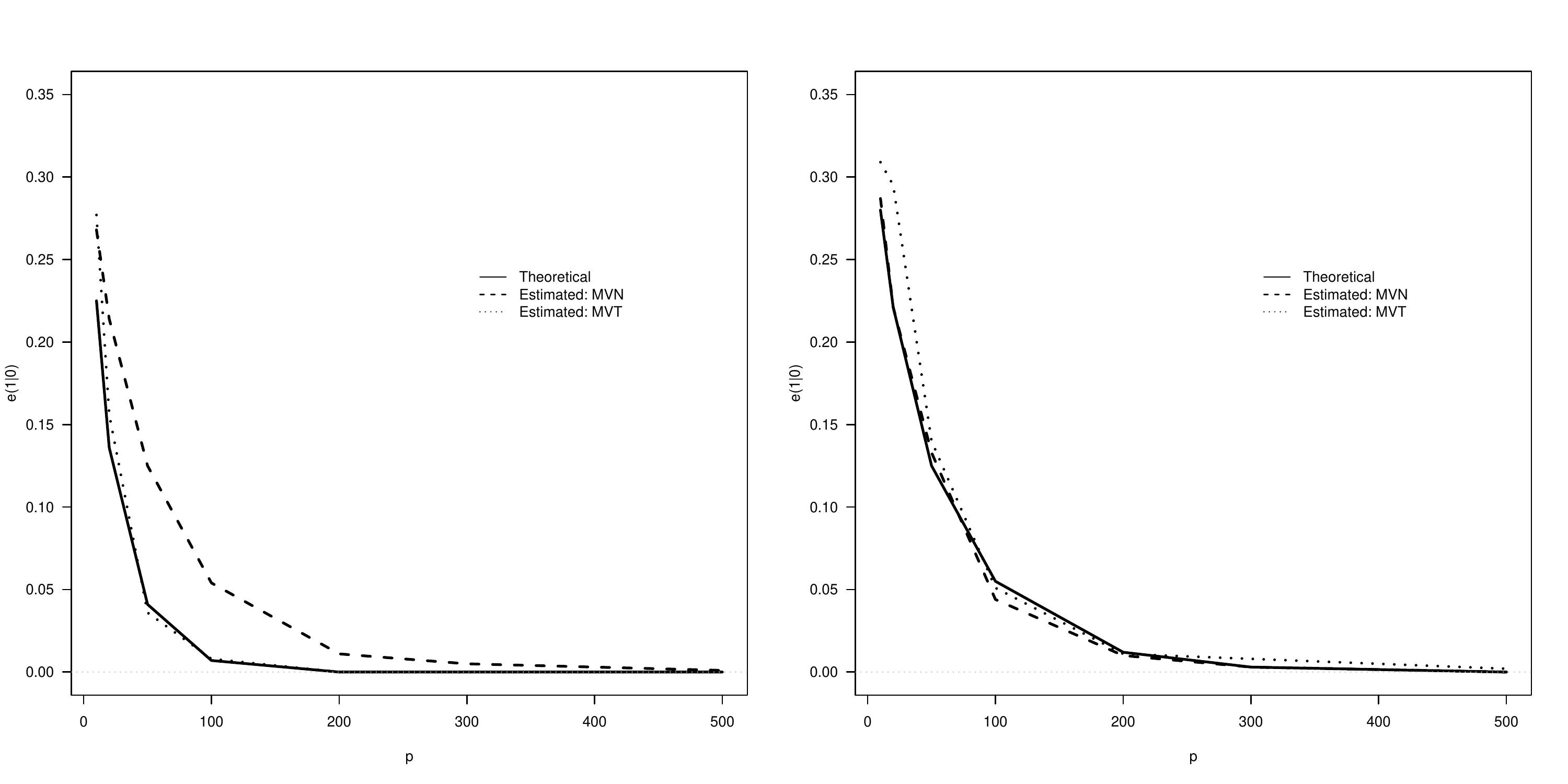}\\
\includegraphics[width = 1.025\textwidth]{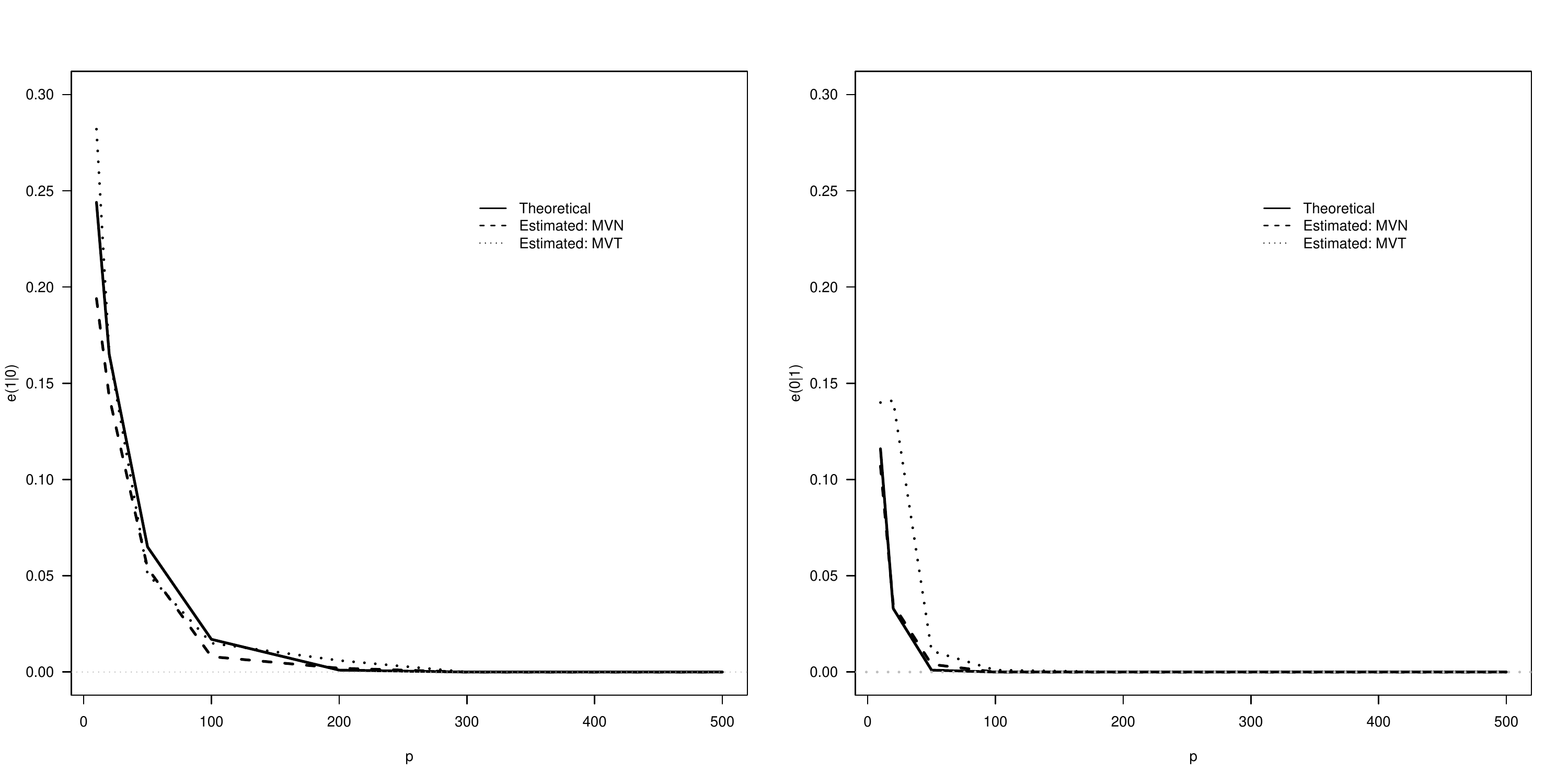}
\caption{\footnotesize \sc Theoretical (thick line) and Estimated error rates of $A_0({\bf x})$ for two-class case with data from multivariate normal (dashed line) and $t$ (dotted line) distributions, $n_1 = 10$, $n_2 = 12$, $p = 10, 20, 50, 100, 200, 300, 500$ and covariance structures AR-AR (upper panel) \& AR-UN (lower penal).}
\end{figure}
distribution with $n_1 = 5, n_2 = 7$, the estimated error rates are relatively higher than under normality, but with $n_i$ increased only by 5, a discernable difference in the performance of the classifier is observed in Fig. 4. Note that, the x-axis in Figs. 3-4 is truncated at $p = 500$ since the misclassification rates already converge to 0 by this value and remain so for larger $p$.

%%%%%%%%%%%%%%%%%%%%%%%%%%%%%%%%%%%%%%%%%%%%%%%%%%%%%%%%%%%%%%%
\section{Applications}\label{sec:Applns}
%%%%%%%%%%%%%%%%%%%%%%%%%%%%%%%%%%%%%%%%%%%%%%%%%%%%%%%%%%%%%%%

We apply $A_0({\bf x})$ on two large data sets for $g = 2$ and $3$. With moderate sample sizes (77 and 102), we use $K = 3$-fold CV for evaluation \citep[see][]{Dudoit02}. \\
%%%%%%%%%%%%%%%%%%%%%%%%%%%%%%%%%%%%%%%%%%%%%%%%
\indent Let $\mathcal{L}$ and $\mathcal{T}$ denote the learning and test sets. We randomly divide data set into $K$ classes of roughly equal size where $\mathcal{T}$ consists of $K - 1$ classes with $K$th class held out as test data. The procedure is repeated $K$ times, each time with a different test class, and a misclassification rate is computed for each repetition. The evaluation criterion is the average misclassification rate over all repetitions. \\
%%%%%%%%%%%%%%%%%%%%%%%%%%%%%%%%%%%%%%%%%%%%%%%%
\indent For $k$th fold of CV, let $n^k_i(\mathcal{L})$, $n^k_i(\mathcal{T})$ and $m^k_{ij}(\mathcal{T})$ be, respectively, the sample sizes for learning and test data in sample $i$ and the number of misclassified observations from class $i$ into class $j$, $i, j = 1, \ldots, g = 2$ or 3, $k = 1, \ldots, K = 3$. Let $e^k(i|j)$ be the estimated misclassification rate, an estimator of $\pi(i|j)$ in (\ref{eqn:MisclRate}), for $k$th rotation, i.e. $e^k(i|j) = m^k_{ij}(\mathcal{T})/n^k(\mathcal{T})$, where $n^k(\mathcal{T}) = n^k_i(\mathcal{T}) + n^k_j(\mathcal{T})$. For $g = 3$, we do the same procedure for each of three pairs and compute overall misclassification rate. For details on the used and other data sets, see \citet{Statnikov05} and also \cite{ShippEtAl02} and \cite{ArmstrongEtAl02}.
%%%%%%%%%%%%%%%%%%%%%%%%%%%%%%%%%%%%%%%%%%%%%%%%
\paragraph{\sc Example 1: DLBCL Data}~The Diffuse Large B-cell Lymphoma (DLBCL) data belongs to a study of lymphoid malignancy in adults. The analysis reported here consists of $p = 5469$ gene expressions studied on pre-treatment biopsies from two independent groups of 77 patients, one with DLBCL ($n_1 = 58$), the other with follicular lymphoma (FL) ($n_2 = 19$). \\
%%%%%%%%%%%%%%%%%%%%%%%%%%%%%%%%%%%%%%%%%%%%%%%%
\indent For a 3-fold CV, we randomly divide the data into three groups of sizes 26, 26, 25 with $n^1(\mathcal{L}) = 52$, $n^1(\mathcal{T}) = 25$ and $n^k(\mathcal{L}) = 51$, $n^k(\mathcal{T}) = 26$ for $k = 2, 3$. By coding the populations as 1 (DLBCL) and 2 (FL), the misclassifications observed from the three rotations of CV, i.e. $m^k_{12}$ and $m^k_{21}$, for $k = 1, 2, 3$, are
\[
m^1_{12} = 3,~ m^1_{21} = 1; \quad m^2_{12} = 6,~ m^2_{21} = 0; \quad m^3_{12} = 2, ~m^3_{21} = 3,
\]
so that the overall misclassification rate is computed as $15/77$. Although of relatively less importance, due to randomly sampled folds for cross-validation, we also report the sample sizes for each fold as following:
\begin{eqnarray*}
  K = 1: && n^1_1(\mathcal{L}) = 38,~n^1_2(\mathcal{L}) = 14,~ n^1_1(\mathcal{T}) = 20,~n^1_2(\mathcal{T}) = 5\\%,~ U^1 = 40999.75 \\
  K = 2: && n^2_1(\mathcal{L}) = 40,~n^2_2(\mathcal{L}) = 11,~ n^2_1(\mathcal{T}) = 18,~n^2_2(\mathcal{T}) = 8\\%,~ U^2 = 187891.3 \\
  K = 3: && n^3_1(\mathcal{L}) = 38,~n^3_2(\mathcal{L}) = 13,~ n^3_1(\mathcal{T}) = 20,~n^3_2(\mathcal{T}) = 6%,~ U^3 = 107175.7.
\end{eqnarray*}
%%%-------------------------------------%%%
\paragraph{\sc Example 2: Leukemia Data}~The data set pertains to a study of patients with acute lymphoblastic leukemia (ALL) carrying a chromosomal translocation involving mixed-lineage leukemia (MLL) gene. The analysis reported here consists of $p = 11225$ gene expression profiles of leukemia cells from $n_2 = 24$ patients diagnosed with B-precursor ALL carrying an MLL translocation and compared to a group of $n_3 = 20$ individual diagnosed with conventional B-precursor without MLL translocation. In addition, there is a third  group of a random sample of $n_1 = 28$ with acute myelogenous leukemia (AML). \\%For details, see \cite{ArmstrongEtAl02}.\\
%%%%%%%%%%%%%%%%%%%%%%%%%%%%%%%%%%%%%%%%%%%%%%%%%
\indent For a 3-fold cross-validation, we randomly divide the data into three equal groups of size 24 each and use $K - 1 = 2$ classes of total $n^k(\mathcal{L}) = 48$ observations in learning set and $n^k(\mathcal{T}) = 24$ in the test set every time, $k = 1, 2, 3$. The rest of the procedure is the same as explained in the Example 1 above. We obtain the following misclassifications for the three folds of cross-validation:
\begin{eqnarray*}
  K = 1: && m^1_{12} = 1,~ m^1_{21} = 0; \quad m^1_{13} = 2,~ m^1_{31} = 1; \quad m^1_{23} = 0, ~m^1_{32} = 1\\
  K = 2: && m^2_{12} = 0,~ m^2_{21} = 1; \quad m^2_{13} = 2,~ m^2_{31} = 0; \quad m^2_{23} = 1, ~m^2_{32} = 0\\
  K = 3: && m^3_{12} = 0,~ m^3_{21} = 0; \quad m^3_{13} = 0,~ m^3_{31} = 0; \quad m^3_{23} = 0, ~m^3_{32} = 0.
\end{eqnarray*}
This gives an overall misclassification rate $9/72$. The sample sizes used in each rotation are also reported below.
\begin{small}
\begin{eqnarray*}
  K = 1:&&n^1_1(\mathcal{L}) = 21,~n^1_2(\mathcal{L}) = 15,~ n^1_3(\mathcal{L}) = 12;~n^1_1(\mathcal{T}) = 7,~n^1_2(\mathcal{T}) = 9,~ n^1_3(\mathcal{T}) = 8\\
  K = 2:&&n^2_1(\mathcal{L}) = 18,~n^2_2(\mathcal{L}) = 14,~ n^2_3(\mathcal{L}) = 16;~n^2_1(\mathcal{T}) = 10,~n^2_2(\mathcal{T}) = 10,~ n^2_3(\mathcal{T}) = 4\\
  K = 3:&&n^3_1(\mathcal{L}) = 17,~n^3_2(\mathcal{L}) = 19,~ n^3_3(\mathcal{L}) = 12;~n^3_1(\mathcal{T}) = 11,~n^3_2(\mathcal{T}) = 5,~ n^3_3(\mathcal{T}) = 8
\end{eqnarray*}
\end{small}

%%%%%%%%%%%%%%%%%%%%%%%%%%%%%%%%%%%%%%%%%%%%%%%%%%%%%%%%%%%%%%%
\section{Discussion and conclusions}\label{sec:Discussion}
%%%%%%%%%%%%%%%%%%%%%%%%%%%%%%%%%%%%%%%%%%%%%%%%%%%%%%%%%%%%%%%

A $U$-classifier for high-dimensional and possibly non-normal data is proposed. The threshold part of the classifier, called $U$-component, is a linear combination of two bivariate $U$-statistics of computed from the two independent samples. The discriminant function part, called $P$-component, forms an inner product between the observation to be classified and the difference of the mean vectors of the corresponding independent samples. It results into a computationally simple classifier which is linear without requiring the underlying covariance matrices to be equal. A multi-class extension with same properties is also given.\\
%%%%%%%%%%%%%%%%%%%%%%%%%%%%%%%%%%%%%%%%%%%%%%%%%%%%%%%%%%%
\indent The classifier is unbiased, consistent and asymptotically normal, under a general multivariate model, including (but not necessarily) the multivariate normal distribution. Rapid convergence of the misclassification rate of the classifier is shown for very small sample sizes and non-normal distributions, under mild and practically justifiable assumptions. The performance of the classifier, in terms of its consistency, asymptotic normality and control of misclassification rate, is shown through simulations for normal and non-normal distributions with sample sizes as small as 5 or 7 and for arbitrary large dimension.\\
%%%%%%%%%%%%%%%%%%%%%%%%%%%%%%%%%%%%%%%%%%%%%%%%%%%%%%%%%%%
\indent We apply the classifier to genetics and microarray data sets, some of the most popular areas for classification analysis. To emphasize the role of high-dimensionality, we demonstrate that the use, accuracy, and validity of the classifier does not rest on any form of data pre-processing as is usually shown in the literature. In other words, a data set measured in large dimension can be directly used for the classifier without any pre-requisites of reducing the dimension through sorting or clustering or other means.
%\newpage
%-----------------------------------
\appendix
%-----------------------------------

%--------------------------------------------------------------------------
\section{Moments of Quadratic and Bilinear Forms}\label{subsec:MomsQFBF}
%--------------------------------------------------------------------------

For ${\bf z}_{ik}$ in (\ref{eqn:ModelLINN}), let $A_{ik} = {\bf z}'_{ik}{\bs \Sigma}_i{\bf z}_{ik} = {\bf y}'_{ik}{\bf A}^2{\bf y}_{ik}$ and $A_{ijkl} = {\bf z}'_{ik}{\bf z}_{jl} = {\bf y}'_{ik}{\bs \Lambda}_i{\bs \Lambda}_j{\bf y}_{jl}$, $k \neq l$, be a quadratic and a bilinear form of independent components with $A_{ik} = {\bf y}'_{ik}{\bf y}_{ik} = Q_{ik}$ for ${\bs \Sigma}_i = {\bf I}$. As all terms involving ${\bf A}_i$ eventually vanish under Assumption \ref{assn:traceSigmaHadProd}, we write ${\bf A}_i = {\bf A}$ for simplicity. Theorem \ref{thm:BasicQFBFMomsNN} gives basic moments of quadratic and bilinear forms which we extend in Lemma \ref{lem:QFBFResults}. Proofs of these results are tedious but simple, therefore skipped \citep[see][]{Ahmad14b}.
%%%%%%%%%%%%%%%%%%%%%%%%%%%%%%%%%%%%%%%%%%%%%%%%%%%%%%%%%%%%%%%
%%%%%%%%%%%%%%%%%%%%%%%%%%%%%%%%%%%%%%%%%%%%%%%%%%%%%%%%%%%%%%%
\begin{thm}\label{thm:BasicQFBFMomsNN} For $A_{ik}$ and $A_{ijkl}$, as defined above, we have
\begin{small}
\begin{eqnarray}
\text{E}\left(Q^2_{ik}\right)^2 &=& 2\text{tr}({\bs \Sigma}^2_i) + [\text{tr}({\bs \Sigma}_i)]^2 + M_1\label{eqn:Mom2QFsingle}\\
\text{E}\left(A^2_{ik}\right)^2 &=& 2\text{tr}({\bs \Sigma}^4_i) + [\text{tr}({\bs \Sigma}^2_i)]^2 + M_2\label{eqn:Mom2QFsingle}\\
\text{E}\left(A_{ik}A_{jk}\right) &=& 2\text{tr}({\bs \Sigma}^3_i{\bs \Sigma}_j) + \text{tr}({\bs \Sigma}^2_i)\text{tr}({\bs \Sigma}_i{\bs \Sigma}_j) + M_2\label{eqn:Mom2QF}\\
\text{E}\left(A^4_{ijkl}\right) &=& 6\text{tr}({\bs \Sigma}_i{\bs \Sigma}_j)^2 + 3\left[\text{tr}({\bs \Sigma}_i{\bs \Sigma}_j)\right]^2 + M_3\label{eqn:Mom4BF}\\
\text{E}\left(Q_{ik}Q_{jk}A^2_{ijkl}\right) &=& 4\text{tr}({\bs \Sigma}_i{\bs \Sigma}_j)^2 + 4\text{tr}({\bs \Sigma}^3_i)\text{tr}({\bs \Sigma}_j) + [\text{tr}({\bs \Sigma}_i)]^2\text{tr}({\bs \Sigma}^2_j) + M_4~~~\label{eqn:QFBFProdMom}
\end{eqnarray}
with $M_1 = \gamma\text{tr}({\bf A}\odot{\bf A})$, $M_2 = \gamma\text{tr}({\bf A}^2\odot{\bf A}^2)$, $M_3 = 6\gamma\text{tr}({\bf A}^2\odot{\bf A}^2) + \gamma^2\sum_{s = 1}^p\sum_{t = 1}^pA^4_{st}$, $M_4 = 2\gamma\text{tr}({\bs \Sigma}_i)\text{tr}({\bf A}^2\odot{\bf A}) + 4\gamma\text{tr}({\bf A}^3\odot{\bf A}) + \gamma\text{tr}({\bf A}\odot{\bf A}{\bf D}{\bf A})$ and ${\bf D} = diag({\bf A})$. Moreover, $\text{E}(A_{ik}) = \text{tr}({\bs \Sigma}_i)$, $\text{E}(A^2_{ikr}) = \text{tr}({\bs \Sigma}^2_i)$ and $\text{Cov}(A_{ik}, A_{ikr}) = 0$. %without normality assumption.
\end{small}
\end{thm}
%%%%%%%%%%%%%%%%%%%%%%%%%%%%%%%%%%%%%%%%%%%%%%%%%%%%%%%%%%%%%%
%%%%%%%%%%%%%%%%%%%%%%%%%%%%%%%%%%%%%%%%%%%%%%%%%%%%%%%%%%%%%%
\begin{lem}\label{lem:QFBFResults} Let ${\bf z}_{it}$ be as given above with ${\bf z}_{it}$, ${\bf z}_{iu}$ independent if $t \neq u$. Then
\begin{eqnarray}
\text{E}[{\bf z}'_{it}{\bf z}_{iu}{\bf z}'_{it}{\bf z}_{iv}{\bf z}'_{iu}{\bs \Sigma}_i{\bf z}_{iv}] &=& \text{tr}({\bs \Sigma}^4_i)\label{eqn:A8}\\
\text{E}[{\bf z}'_{it}{\bf z}_{iu}{\bf z}'_{iw}{\bf z}_{iu}{\bf z}'_{it}{\bf z}_{iv}{\bf z}'_{iw}{\bf z}_{iv}] &=& \text{tr}({\bs \Sigma}^4_i)\label{eqn:A9}\\
\text{E}[({\bf z}'_{it}{\bf z}_{iu})^2{\bf z}'_{it}{\bs \Sigma}_i{\bf z}_{it}] &=& 2\text{tr}({\bs \Sigma}^4_i) + [\text{tr}({\bs \Sigma}^2_i)]^2 + M_2\label{eqn:A10}\\
\text{Var}({\bf z}'_{it}{\bf z}_{iu}{\bf z}'_{iv}{\bf z}_{iu}) &=& 2\text{tr}({\bs \Sigma}^4_i) + [\text{tr}({\bs \Sigma}^2_i)]^2 + M_2\label{eqn:A11}\\
\text{Cov}[({\bf z}'_{it}{\bf z}_{iu})^2, ({\bf z}'_{it}{\bf z}_{iv})^2] &=& 2\text{tr}({\bs \Sigma}^4_i) + M_2\label{eqn:A12}\\
\text{E}[({\bf z}'_{it}{\bf z}_{ju})^2{\bf z}'_{it}{\bs \Sigma}_j{\bf z}_{it}] &=& 2\text{tr}\{({\bs \Sigma}_i{\bs \Sigma}_j)^2\} + \big[\text{tr}({\bs \Sigma}_i{\bs \Sigma}_j)\big]^2 + M_2~~\label{eqn:A13}\\
\text{Var}({\bf z}'_{it}{\bf z}_{ju}{\bf z}'_{iv}{\bf z}_{ju}) &=& 2\text{tr}\{({\bs \Sigma}_i{\bs \Sigma}_j)^2\} + \big[\text{tr}({\bs \Sigma}_i{\bs \Sigma}_j)\big]^2 + M_2~~\label{eqn:A14}\\
\text{Cov}[({\bf z}'_{jt}{\bf z}_{iu})^2, ({\bf z}'_{jt}{\bf z}_{iv})^2] &=& 2\text{tr}\{({\bs \Sigma}_i{\bs \Sigma}_j)^2\} + M_2\label{eqn:A15}\\
\text{E}({\bf z}'_{jt}{\bf z}_{iu}{\bf z}'_{jt}{\bf z}_{iv}{\bf z}'_{iu}{\bs \Sigma}_j{\bf z}_{iv}) &=& \text{tr}\{({\bs \Sigma}_i{\bs \Sigma}_j)^2\}\label{eqn:A16}\\
\text{Cov}[({\bf z}'_{it}{\bf z}_{iu})^2, {\bf z}'_{it}{\bs \Sigma}_j{\bf z}_{it}] &=& 2\text{tr}({\bs \Sigma}^3_i{\bs \Sigma}_j) + M_2\label{eqn:A17}\\
\text{Cov}({\bf z}'_{it}{\bs \Sigma}_i{\bf z}_{iu}, {\bf z}'_{it}{\bs \Sigma}_j{\bf z}_{iu}) &=& \text{tr}\{({\bs \Sigma}_i{\bs \Sigma}_j)^2\}\label{eqn:A18}\\
\text{E}[({\bf z}'_{iu}{\bf z}_{iv})^2{\bf z}'_{it}{\bs \Sigma}_j{\bf z}_{it}] &=& \text{tr}({\bs \Sigma}_i{\bs \Sigma}_j)\text{tr}({\bs \Sigma}_i)^2,\label{eqn:A19}
\end{eqnarray}
where $\text{E}[({\bf z}'_{it}{\bf z}_{iu})^2{\bf z}'_{it}{\bf z}_{iu}{\bf z}'_{it}{\bf z}_{iv}]$, $\text{E}[{\bf z}'_{it}{\bf z}_{iu}{\bf z}'_{it}{\bf z}_{iv}{\bf z}'_{it}{\bs \Sigma}_i{\bf z}_{it}]$, $\text{E}[{\bf z}'_{it}{\bf z}_{iu}{\bf z}'_{it}{\bf z}_{iv}{\bf z}'_{it}{\bs \Sigma}_i{\bf z}_{iu}]$,\\ $\text{E}[({\bf z}'_{it}{\bf z}_{iu})^2{\bf z}'_{it}{\bs \Sigma}_i{\bf z}_{iu}]$, $\text{E}[({\bf z}'_{it}{\bf z}_{iu})^2{\bf z}'_{it}{\bs \Sigma}_j{\bf z}_{iu}]$, $\text{E}[({\bf z}'_{it}{\bf z}_{iu})^2{\bf z}'_{it}{\bf z}_{iv}{\bf z}'_{iu}{\bf z}_{iv}]$ all vanish.
\end{lem}
%%%%%%%%%%%%%%%%%%%%%%%%%%%%%%%%%%%%%%%%%%%%%%%%%%%%%%%%%%%%%%%
%%%%%%%%%%%%%%%%%%%%%%%%%%%%%%%%%%%%%%%%%%%%%%%%%%%%%%%%%%%%%%%

%%%%%%%%%%%%%%%%%%%%%%%%%%%%%%%%%%%%%%%%%%%%%%%%%%%%%%%
\section{Main Proofs}\label{sec:MainProofs}
%%%%%%%%%%%%%%%%%%%%%%%%%%%%%%%%%%%%%%%%%%%%%%%%%%%%%%%

%---------------------------------------------------------------------------------------%
\subsection{Proof of Lemma \ref{lem:MomsTwoSDF}}\label{subsec:ProofLemmaMomsTwoSDF}
%---------------------------------------------------------------------------------------%

Let ${\bf x}_0 \in \pi_1$. With ${\bf x}_0$ independent of both samples, $\text{E}[A_0({\bf x})]$ is trivial. For variance, ignoring $p$ for simplicity, we begin with
\[
\text{Var}[{\bf x}'_0(\overline{\bf x}_1 - \overline{\bf x}_2)] = \text{E}[{\bf x}'_0(\overline{\bf x}_1 - \overline{\bf x}_2)]^2 - [{\bs \mu}'_1({\bs \mu}_1 - {\bs \mu}_2)]^2.
\]
Since $\text{E}[{\bf x}'_0(\overline{\bf x}_1 - \overline{\bf x}_2)]^2 = \text{E}[{\bf x}'_0(\overline{\bf x}_1 - \overline{\bf x}_2)(\overline{\bf x}_1 - \overline{\bf x}_2)'{\bf x}_0]$, we immediately get
\begin{eqnarray*}
\text{E}[{\bf x}'_0(\overline{\bf x}_1 - \overline{\bf x}_2)]^2 &=& \text{tr}\left[({\bs \Sigma}_1 + {\bs \mu}_1{\bs \mu}'_1)\left(\frac{{\bs \Sigma}_1}{n_1} + \frac{{\bs \Sigma}_2}{n_2} + ({\bs \mu}_1 - {\bs \mu}_2)({\bs \mu}_1 - {\bs \mu}_2)'\right)\right],
\end{eqnarray*}
so that
\[
\text{Var}[{\bf x}'_0(\overline{\bf x}_1 - \overline{\bf x}_2)] = \frac{\text{tr}({\bs \Sigma}^2_1)}{n_1} + \frac{\text{tr}({\bs \Sigma}_1{\bs \Sigma}_2)}{n_2} + \frac{{\bs \mu}'_1{\bs \Sigma}_1{\bs \mu}_1}{n_1} + \frac{{\bs \mu}'_1{\bs \Sigma}_2{\bs \mu}_1}{n_2} + ({\bs \mu}_1 - {\bs \mu}_2)'{\bs \Sigma}_1({\bs \mu}_1 - {\bs \mu}_2).
\]
Now $\text{Var}(U_{n_1} - U_{n_2}) = \sum_{i = 1}^2\text{Var}(U_{n_i})$. For $U_{n_i}$ with $h({\bf x}_{ik}, {\bf x}_{ir}) = {\bf x}'_{ik}{\bf x}_{ir}$ \citep[][Ch. 5]{Serfling80}), $h_1({\bf x}_{ik}) = {\bf x}'_{ik}{\bs \mu}_i$ with $\xi_1 = \text{Var}[h_1({\bf x}_{ik})] = {\bs \mu}'_i{\bs \Sigma}_i{\bs \mu}_i$, and $h_2(\cdot) = h(\cdot)$ with $\xi_2 = \text{Var}(A_{ik}) = \text{tr}({\bs \Sigma}^2_i) + 2{\bs \mu}'_i{\bs \Sigma}_i{\bs \mu}_i$, so that
\[
\text{Var}(U_{n_i}) = \frac{2}{n_i(n_i - 1)}[2(n_i - 2)\xi_1 + \xi_2] = \frac{2\text{tr}({\bs \Sigma}^2_i)}{n_i(n_i - 1)} + \frac{4{\bs \mu}'_i{\bs \Sigma}_i{\bs \mu}_i}{n_i}, ~i = 1, 2.
\]
For $\text{Cov}[{\bf x}'_0(\overline{\bf x}_1 - \overline{\bf x}_2), U_{n_1} - U_{n_2}]$, $\text{Cov}({\bf x}'_0\overline{\bf x}_2, U_{n_1}) = 0 = \text{Cov}({\bf x}'_0\overline{\bf x}_1, U_{n_2})$, by independence, where it immediately follows that $\text{Cov}({\bf x}'_0\overline{\bf x}_i, U_{n_i}) = 2{\bs \mu}'_i{\bs \Sigma}_i{\bs \mu}_i/n_i$, $i = 1, 2$. Combining all results and simplifying gives Eqn. (\ref{eqn:VarDFTwoS}).

%---------------------------------------------------------------------------------------%
\subsection{Proof of Theorem \ref{thm:PropsEstrsBounds}}\label{subsec:ProofThmPropsEstrs}
%---------------------------------------------------------------------------------------%

The unbiasedness is trivial. For $\text{Var}(E_0)$, $\text{Var}(U_{n_i})$, $i = 1, 2$, are given in Sec. \ref{subsec:ProofLemmaMomsTwoSDF}. For $\text{Var}(U_{n_{ij}})$, $h({\bf x}_{ik}, {\bf x}_{jl}) = {\bf x}'_{ik}{\bf x}_{jl}$ with $h_{10} = {\bs \mu}'_j{\bf x}_{ik}$, $h_{01} = {\bs \mu}'_i{\bf x}_{jl}$ so that $\xi_{10} = \text{Var}[h_{10}(\cdot)] = {\bs \mu}'_j{\bs \Sigma}_i{\bs \mu}_j$ and $\xi_{10} = \text{Var}[h_{10}(\cdot)] = {\bs \mu}'_i{\bs \Sigma}_j{\bs \mu}_i$. Also $h_{11}(\cdot) = h(\cdot)$ with $\xi_{11} = \text{Var}[h_{11}(\cdot)] = {\bs \mu}'_i{\bs \Sigma}_j{\bs \mu}_i + {\bs \mu}'_j{\bs \Sigma}_i{\bs \mu}_j + \text{tr}({\bs \Sigma}_i{\bs \Sigma}_j)$. Hence \citep{Lee90}
\[
\text{Var}(U_{n_in_j}) = \frac{1}{n_in_jp^2}\left[n_i{\bs \mu}'_i{\bs \Sigma}_j{\bs \mu}_i + n_j{\bs \mu}'_j{\bs \Sigma}_i{\bs \mu}_j + \text{tr}({\bs \Sigma}_i{\bs \Sigma}_j)\right]
\]
where $\text{Cov}(U_{n_i}, U_{n_in_j}) = 2{\bs \mu}'_j{\bs \Sigma}_i{\bs \mu}_i/n_ip^2$, $\text{Cov}(U_{n_j}, U_{n_in_j}) = 2{\bs \mu}'_i{\bs \Sigma}_j{\bs \mu}_j/n_jp^2$ and $\text{Cov}(U_{n_i}, U_{n_j}) = 0$ by independence. $\text{Var}(E_0/p)$ can now be approximated as
\[
\text{Var}(E_0/p) = 2\text{tr}({\bs \Sigma}^2_{0ij})/p^2 + 4({\bs \mu}_i - {\bs \mu}_j)'{\bs \Sigma}_{0ij}({\bs \mu}_i - {\bs \mu}_j)/p^2,~ i, j = 1, 2,~i \neq j.
\]
With second term vanishing under Assumption \ref{assn:ExtraH1Distn} and first term bounded in $p$ under Assumption \ref{assn:traceSigma}, $\text{Var}(E_0)$ reduces to $O(1/n_1 + 1/n_2)^2$ as $p \rightarrow \infty$, so that the consistency follows immediately as $n_i, p \rightarrow \infty$. The bound in (\ref{eqn:BoundVarE0}) also follows trivially. As $E_i$ and $E_{ij}$, are also one- and two-sample $U$-statistics with higher order kernels, we essentially follow the same strategy as for $E_0$. First, from Theorem \ref{thm:BasicQFBFMomsNN} and Lemma \ref{lem:QFBFResults}, it can be shown that \citep[see][Ch. 2]{Ahmad14b}
\begin{eqnarray}
\text{Var}(E_i) &=& \frac{4}{\eta(n_i)p^4}\Big[(2n^3_i - 9n^2_i + 9n_i - 16)\text{tr}({\bs \Sigma}^4_i) \nonumber\\
&&+~(n^2_i - 3n_i + 8)[\text{tr}({\bs \Sigma}^2_i)]^2 + M_2O(n^3_i) + M_3O(n^2_i)\Big]\\
\text{Var}(E_{ij}) &=& \frac{2}{(n_i - 1)(n_j - 1)p^4}\Big[(n - 1)\text{tr}\{({\bs \Sigma}_1{\bs \Sigma}_2)^2\} + [\text{tr}({\bs \Sigma}_1{\bs \Sigma}_2)]^2 \nonumber\\
&&\qquad\qquad\qquad\qquad\qquad\qquad+~M_2O(n) + M_3O(1)\Big]\\
\text{Cov}(E_i, E_{ij}) &=& \frac{4}{n_i(n_i - 1)p^4}\Big[n_i\text{tr}({\bs \Sigma}^3_i{\bs \Sigma}_j) + M_2O(n_i)\Big]
\end{eqnarray}
$n = n_i + n_j$, $i, j = 1, 2$, $i \neq j$, $M_2$, $M_3$ are as in Theorem \ref{thm:BasicQFBFMomsNN} and $\text{Cov}(E_i, E_j) = 0$. As terms involving $M$'s vanish under Assumption \ref{assn:traceSigmaHadProd}, the consistency and the bounds (by Cauchy-Schwarz inequality) follow the same way as for $E_0$. Note also that, the terms involving $M$'s are exactly zero under normality in which case the same results follow even more conveniently.

%---------------------------------------------------------------------------------------%
\subsection{Proof of Theorem \ref{thm:ConsistencyTwoSDF}}\label{subsec:ProofThmConsisTWoSDF}
%---------------------------------------------------------------------------------------%

The proof essentially follows from that of Theorem \ref{thm:PropsEstrsBounds} without much new computations. In particular, the first part, assuming true parameters known, is trivial. For the second part with empirical estimators, the $(n_i, p)$-consistency of estimators proved in Sec. \ref{subsec:ProofThmPropsEstrs} implies that $E_0/\text{E}(E_0) \xrightarrow{\mathcal{P}} 1$, and the same holds for $E_i$, $E_{ij}$. Plugging these estimators in the moments of $A_0({\bf x})$ and using Slutsky's lemma, $\widehat{\delta}^2_i/\delta^2_i \xrightarrow{\mathcal{P}} 1$ so that $\widehat{\text{Var}[A_0({\bf x})]} = \text{Var}[A_0({\bf x}] + o_P(1)$, and the consistency follows similarly as with known parameters.

%---------------------------------------------------------------------------------------%
\subsection{Proof of Theorem \ref{thm:AsympNTwoSDF}}\label{subsec:ProofThmAsympNTwoSDF}
%---------------------------------------------------------------------------------------%

Write $\widetilde{A}_0({\bf x}) = [A_0({\bf x})|{\bf x}_0 \in \pi_1]- \text{E}[A_0({\bf x})|{\bf x}_0 \in \pi_1]$, where
\[
\widetilde{A}_0({\bf x}) = [{\bf x}'_0(\overline{\bf x}_1 - \overline{\bf x}_2) - {\bs \mu}'_1({\bs \mu}_1 - {\bs \mu}_2)] - [(U_{n_1} - {\bs \mu}'_1{\bs \mu}_1) - (U_{n_2} - {\bs \mu}'_2{\bs \mu}_2)]/2,
\]
ignoring $p$ for simplicity. Let $\widehat{U}_{n_i}$ be the projection of $\widetilde{U}_{n_i} = U_{n_i} - {\bs \mu}'_i{\bs \mu}_i$, $i = 1, 2$. Then \citep[][Ch. 5]{Serfling80} $g_1({\bf x}_{1k}) = h_1(\cdot) - {\bs \mu}'_1{\bs \mu}_1 = ({\bf X}_{1k} - {\bs \mu}_1)'{\bs \mu}_1$ for $U_{n_1}$, and similarly $g_1({\bf x}_{2l})$ for $U_{n_2}$, with $\text{E}[g(\cdot)] = 0$ in both cases, so that
\[
\widehat{U}_{n_1} - \widehat{U}_{n_2} = \frac{2}{n_1}\sum_{k = 1}^{n_1}({\bf X}_{1k} - {\bs \mu}_1)'{\bs \mu}_1 - \frac{2}{n_2}\sum_{l = 1}^{n_2}({\bf X}_{2l} - {\bs \mu}_2)'{\bs \mu}_2,
\]
where $\widetilde{U}_{n_i} = \widehat{U}_{n_i} + o_P(1)$, $i = 1, 2$. With ${\bf x}_0 \in \pi_1$ and independence of samples, this projection of $\widetilde{A}_0({\bf x})$ results into a sum of two independent components, each an average of iid variables \citep{vdv98}. Taking $p$ into account, the asymptotic normality follows by the CLT under Assumptions \ref{assn:4thMomnt}-\ref{assn:ExtraH1Distn} as $n_i, p \rightarrow \infty$.

%%%%%%%%%%%%%%%%%%%%%%%%%%%%%%%%%%%%%%%%%%%%%%%%
\paragraph{Acknowledgement} We are sincerely thankful to {\sc Prof. Thomas} {\sc Mikosch}, Department of Mathematics, University of Copenhagen, for careful perusal of the paper, discussion and encouragement. The research of Tatjana Pavlenko is partially supported by the Swedish Research Council, Grant No. C0595201.

%%%%%%%%%%%%%%%%%%%%%%%%%%%%%%%%%%%%%%%%%%%%%%%%%%%%%%%%%%%%%%%%%%%%%%%%%%%%%%%%%%%%%%%%%%%%%%%%%%%%%%%%%%%%%%%%%%%%%%%%%%%%%%%%%%%%%%%%%%%%%%%%%
%%%%%%%%%%%%%%%%%%%%%%%%%%%%%%%%%%%%%%%%%%%%%%%%%%%%%%%%%%%%%%%%%%%%%%%%%%%%%%%%%%%%%%%%%%%%%%%%%%%%%%%%%%%%%%%%%%%%%%%%%%%%%%%%%%%%%%%%%%%%%%%%%

\end{document}